\documentclass[12pt]{article}

    \newcommand{\lab}[1]{\label{#1}}                

\usepackage{amssymb}
\usepackage{amsfonts}
\usepackage{mathrsfs}
\usepackage{graphicx}
\usepackage{amsmath}
\usepackage{epsfig}
\usepackage{verbatim}

\def\smallpage{
\addtolength\textwidth{3cm} \addtolength\oddsidemargin{-1.5cm}
\addtolength\textheight{3cm} \addtolength\topmargin{-1.5cm}}
\smallpage

\newcommand{\remove}[1]{}
\newcommand{\bel}[1]{\be\lab{#1}}
\newcommand\eqn[1]{(\ref{#1})}

\newcommand{\be}{\begin{equation}}
\newcommand{\ee}{\end{equation}}
\newcommand{\bea}{\begin{eqnarray}}
\newcommand{\eea}{\end{eqnarray}}
\newcommand{\bean}{\begin{eqnarray*}}
\newcommand{\eean}{\end{eqnarray*}}
\newtheorem{thm}{Theorem}[section]
\newtheorem{cor}[thm]{Corollary}

\newtheorem{lemma}[thm]{Lemma}

\def\proof{\noindent{\bf Proof.\ }  }
\def\qed{~~\vrule height8pt width4pt depth0pt}

\def\ex{{\bf E}}
\def\pr{{\bf P}}

\def \PA{\mathcal{P}}

\def\mathC{{\mathcal C}}

\def\d{d_{\max}}

\def\Gnp{{\mathcal G}(n,p)}

\def\mathscrC{U}

\def\eps{\epsilon}

\def\ss{\smallskip}
\def\non{\nonumber}
\def\no{\noindent}

\catcode`@=11 \@addtoreset{equation}{section}

\catcode`@=12

\title{Induced subgraphs in sparse random graphs with given degree sequence}
\author{Pu Gao\thanks{Research supported by the Humboldt Foundation}
\\ Max-Planck-Institut f\"{u}r Informatik\\ janegao@mpi-inf.mpg.de \and  Yi Su\footnote{Current affiliation: University of Michigan, MI, USA}\\
University of Waterloo\\ yisu@umich.edu \and Nicholas
Wormald\thanks{Research
supported by the Canadian Research Chairs Program and NSERC}\\
 University of
Waterloo\\nwormald@uwaterloo.ca}
\date{}

\begin{document}
\maketitle

\begin{abstract}
Let $\mathcal{G}_{n,d}$ denote the uniformly random $d$-regular graph on $n$ vertices. For any $S\subset [n]$, we obtain estimates of the probability that the subgraph
of $\mathcal{G}_{n,d}$  induced by $S$ is a given graph $H$. The estimate gives an asymptotic formula for any $d=o(n^{1/3})$, provided that $H$ does not contain almost all the edges of the random graph. The result is
further extended to  the probability space
of random graphs with a given degree sequence.
\end{abstract}

\section{Introduction}

Properties of subgraphs and induced subgraphs in random graph models
have been investigated by various authors. Ruci{\'n}ski~\cite{R,R3}
studied the distribution of the count of small subgraphs in the
standard random graph model $\mathcal{G}_{n,p}$, and  conditions
under which the distribution converges to the normal distribution.
He also studied properties of induced subgraphs in~\cite{R2}.

Techniques for analysing  the standard random graph
model $\mathcal{G}_{n,p}$ often do not apply in the random regular graph
model $\mathcal{G}_{n,d}$. We take the vertex set of the graph to be $[n]$ in both these models. For $S\subseteq [n]$, let $G_S$ denote the subgraph of $G$ induced by $S$. For a graph $H$ with   vertex set $S$, computing the probabilities $\pr(G_S\supseteq H)$ and $\pr(G_S=H)$  in
$\mathcal{G}_{n,p}$ is trivial, but computing them in
$\mathcal{G}_{n,d}$ is not easy, especially when the degree $d\to\infty$ as $n\to\infty$. McKay~\cite{M3} estimated lower
and upper bounds of $\pr(G_S\supseteq H)$ in $\mathcal{G}_{n,d}$
when the degree sequence of $H$ and $d$ satisfy certain conditions.
These bounds are useful in estimating the asymptotic value of
$\pr(G_S\supseteq H)$ when $d$ is not too large or $H$ is small.
Z. Gao and the third author~\cite{GW3} proved that
the distribution of the number of small subgraphs with certain
restrictions (such as $d$ not growing too quickly) converges to the normal distribution in
$\mathcal{G}_{n,d}$. No such results on induced subgraphs
have been derived, although  the main results of~\cite{M3} could be used as a basis for obtaining results on induced subgraphs. However, this would require severe restrictions on the size of the subgraphs, and seems unlikely to apply to subgraphs with more than $n^{2/3}$ vertices for any $d$.

On the other hand, for very dense regular graphs,
Krivelevich, Sudakov and Wormald~\cite{KSW} computed $\pr(G_S= H)$
in $\mathcal{G}_{n,d}$ when $n$ is odd, $d=(n-1)/2$ and
$|V(H)|=o(\sqrt{n})$. McKay~\cite{M4} has recently given a stronger result, for more general degree sequences and provided $H$ has less than $n^{1+\eps}$ edges for some $\eps>0$.

An asymptotic formula of the probability that $G_S = H$ or
$G_S\supseteq H$ in a random bipartite graph with a specified degree
sequence has been derived by Bender~\cite{B3} when the maximum
degree is bounded. The result was extended further by Bollob\'{a}s
and McKay~\cite{BM} and by McKay~\cite{M2} when the maximum degree
goes to infinity slowly as $n$ goes to infinity. Greenhill and
Mckay~\cite{GM} recently derived an asymptotic formula for the case
when the random bipartite graph is sufficiently dense  and $H$ is
sparse enough.

For a vector ${\bf d}=(d_1,\ldots,d_n)$  of nonnegative integers,
let $M=M({\bf d})=\sum_{i=1}^n d_i$ and let $\mathcal{G}_{ {\bf
d}}$ denote the class of graphs with degree sequence ${\bf d}$ and
the uniform distribution (so $\mathcal{G}_{ {\bf d}}$ is a
generalisation of $\mathcal{G}_{n,d}$). In this paper, we compute
the probability that $G_S = H$ in $\mathcal{G}_{ {\bf d}}$ when
$\d=o((M-2m(H))^{1/4})$, where $m(H)$ denotes the number of edges in $H$ and $\d=\max\{d_1,\ldots, d_n\}$.  The power of
this result is that there is no major restriction on the size or density
of $H$.  In
Section~\ref{s:results}, as a direct application of our main result,
we compute the probability that a given set of vertices in
$\mathcal{G}_{n,d}$ is an independent set. Our results will also be useful as a basic tool
for studying the properties of induced subgraphs in the binomial random graph $\Gnp$, such as the subgraph induced by the vertices of even degree, or odd degree.

A graph $G$ is called a {\em B-graph with vertex bipartition $(L,R)$}  if $V(G)=L\cup R$, and $L$ is an independent set of $G$.  If the graph is not necessarily simple, i.e.\ loops and multiple edges are allowed, we call it a {\em B-multigraph} instead.  An edge in a B-graph or B-multigraph is called a {\em
mixed edge} if its end vertices are in $L$ and $R$ respectively, and  a {\em pure edge} if they are both in
$R$. Given a nonnegative integer vector ${\bf d}$, let
$\mathcal{G}(L,R,{\bf d})$ be the set of B-graphs with bipartition
$L$ and $R$ and the degree sequence ${\bf d}$ and let $g(L,R,{\bf
d})=|\mathcal{G}(L,R,{\bf d})|$. By convention,  $g(L,R,{\bf
d})=0$ if ${\bf d}$ is not nonnegative.

Given a sequence ${\bf d}$, let $g({\bf d})$ denote the number of
graphs on vertex set $[n]$ with degree sequence ${\bf d}$. Given
$S=[s]\subset [n]$, let $H$ be a given graph on vertex set $S$ with
degree sequence $(k_i)_{1\le i\le s}$. Let $\bf{d'}$ be the integer
vector defined by $d'_i=d_i-k_i$ for $i\in S$ and $d'_i=d_i$ for
$i\in [n]\setminus S$. Then the number of graphs with degree
sequence ${\bf d}$ and with $G_S = H$ is $g(S,[n]\setminus S,{\bf
d'})$, and so the probability that $G_S = H$ in $\mathcal{G}_{n,{\bf
d}}$ equals $g(S,[n]\setminus S,{\bf d'})/g({\bf d})$. So the study
of induced subgraphs leads directly to the question of counting
B-graphs.

The following theorem by McKay~\cite{M2} gives an asymptotic formula
for  $g({\bf d})$ when $\d^4=o(M({\bf d}))$. (The restriction on
$\d$ was relaxed further by McKay and Wormald in~\cite{MW2}, but to
do so requires a few extra terms in the exponential factor of the
asymptotic formula, and is not needed for the purpose of this
paper.)
\begin{thm}[McKay]\lab{t:regularCount}
Let  ${\bf d}=(d_1,\ldots, d_n)$ with  $\sum_{i=1}^nd_i$ even and
$\d=o(M({\bf d})^{1/4})$. The number of graphs with degree sequence
$\bf d$ is uniformly
$$
\frac{M({\bf d})!}{2^{M({\bf d})/2}(M({\bf d})/2)!\prod_{i=1}^n
d_i!}\cdot\exp\left(-\mu({\bf d})-\mu({\bf d})^2+O(\d^4/M({\bf
d}))\right)
$$
as $n\to\infty$.
\end{thm}
 By ``uniformly'' in the above theorem we mean the constant
implicit in $O(.)$ is the same for all choices of ${\bf d}$ as a
function of $n$, for a given function implicit in the $o(.)$ term. A
special case of Theorem~\ref{t:regularCount} gives that the
 number of $d$-regular graphs on $n$ vertices is asymptotically
$$
\frac{(dn)!}{2^{dn/2}(dn/2)!(d!)^n}\cdot\exp\left(-\frac{d^2-1}{4}\right),
$$
when $d=o(n^{1/3})$.

Our main result is an
asymptotic formula for  $g(L,R,{\bf
d})$, to an accuracy matching McKay's formula in Theorem~\ref{t:regularCount}. This is given in Section~\ref{s:results},
together with its direct applications to estimating $\pr(G_S = H)$ in $\mathcal{G}_{ {\bf d}}$, and some special cases are also given there. The proofs use the
switching method, first introduced by McKay~\cite{M2}, with
refinements by McKay and Wormald~\cite{MW2}, and suitably modified
for our purposes here.  In Section~\ref{s:basics} we use
switchings to estimate the ratios between probabilities defined by
the counts of loops and various types of multiple edges. In
Section~\ref{s:moreswitchings} we again use switchings to evaluate
some variables appearing in those estimates, and in
Section~\ref{s:synthesis} we use these to prove the main theorem.

\section{Main results}
\lab{s:results}

Our main goal in this paper is to estimate $g(L,R,{\bf d})$. We
first define some notation. For any positive integer $n$, let $[n]$
denote the set $\{1,2,\ldots, n\}$. Given a sequence ${\bf
d}=(d_1,\ldots,d_n)$, let $\d=\max\{d_i, i\in[n]\}$ and let
$M_2({\bf d})=\sum_{i=1}^n d_i(d_i-1)$. Define $\mu({\bf d})$ to be
$M_2({\bf d})/2M({\bf d})$.

For any $S\subset L\cup R$, define
$$
M_1({\bf d},S)=\sum_{i\in S}d_i,\quad
M_2({\bf d},S)=\sum_{i\in S}d_i(d_i-1),
$$
\begin{eqnarray}
\mu_0({\bf d},L,R)&=&\frac{(M_1({\bf d},R)-M_1({\bf d},L))M_2({\bf d},R)}{2M_1({\bf d},R)^2},\lab{eq:mu0}\\
\mu_1({\bf d},L,R)&=&\frac{M_2({\bf d},R)M_2({\bf d},L)}{2M_1({\bf
d},R)^2},\lab{eq:mu1}\\
 \mu_2({\bf d},L,R)&=&\mu_0({\bf d},L,R)^2.\lab{eq:mu2}
\end{eqnarray}
We drop the notations $L$ and $R$ from $\mu_i({\bf
d},L,R)$ for $i=0,1,2$ when the context is clear. Note also that if
$M_1({\bf d},R)<M_1({\bf d},L)$, then $g(L,R,{\bf d})$ is trivially
$0$, so we may assume that
\begin{equation}
M_1({\bf d},R)\ge M_1({\bf d},L). \lab{eq:assumption}
\end{equation}

The following theorem, proved in Section~\ref{s:synthesis}, gives an asymptotic formula for
$g(L,R,\bf{d})$.
\begin{thm}\lab{t:graphCount} Let  ${\bf d}=(d_1,\ldots, d_n)$ with   $\sum_{i=1}^n d_i$ even, $\d=o(M({\bf d})^{1/4})$ and $M_1({\bf d},R)\ge M_1({\bf
d},L)$. Then uniformly over all $L$ and ${\bf d}$ as $n\to\infty$,
\begin{eqnarray*}
g(L,R,{\bf d})&=&\frac{M_1({\bf d},R)!e^{-\mu_0({\bf
d})-\mu_1({\bf d})-\mu_2({\bf d})}}{2^{(M_1({\bf
d},R)-M_1({\bf d},L))/2}((M_1({\bf d},R)-M_1({\bf
d},L))/2)!\prod_{i=1}^nd_i!} \left(1 +O\left(\frac{\d^4}{M({\bf d})}\right)\right).
\end{eqnarray*}

\end{thm}
 Applying
Theorems~\ref{t:graphCount} and~\ref{t:regularCount} we directly get
the following. Here $d'_{\max}$ denotes $\max\{d_1',\ldots, d_n'\}$.
\begin{cor}\lab{t2:graphCount}
Let  ${\bf d}=(d_1,\ldots, d_n)$ with   $\sum_{i=1}^n d_i$ even and
$\d=o(M({\bf d})^{1/4})$. Let $S=[s]\subset [n]$, let $H$ be a graph
on vertex set $S$ with degree sequence ${\bf k}=(k_1,\ldots,k_s)$,
let $h=\sum_{i=1}^h k_i$ and let ${\bf d}'=(d'_1,\ldots,d'_n)$ with
$d'_i=d_i-k_i$ for $i\in S$ and $d'_i=d_i$ for $i\notin S$.  If
$d'_i<0$ for some $i\in [n]$ or $M_1({\bf d}',[n]\setminus S)<M_1({\bf
d}',S)$, then $\pr_{\mathcal{G}_{ {\bf d}}}(S,H)=0$. Otherwise, if
$d'_{\max}=o(M({\bf d'})^{1/4})$, then uniformly
\begin{eqnarray*}
\pr_{\mathcal{G}_{ {\bf d}}}(S,H) &= &\exp\left(-\mu_0({\bf
d}')-\mu_1({\bf d}')-\mu_2({\bf d}')+\mu({\bf d})+\mu({\bf
d})^2+O\left(\frac{{d'}_{\max}^4}{M({\bf d'})}+\frac{\d^4}{M({\bf
d})}\right)\right)\\
 &&\times \prod_{i=1}^s[d_i]_{k_i}\frac{M_1({\bf
d}',[n]\setminus S)!2^{M_1({\bf d}',S)+h/2}(M({\bf d})/2)!}{((M_1({\bf
d}',[n] \setminus  S)-M_1({\bf d}',S))/2)!M({\bf d})!}.
\end{eqnarray*}
where $\mu_i({\bf d}')= \mu_i({\bf d}',S,[n]\setminus S)$ for $i=0$, $1$ and $2$.
\end{cor}
\proof Recall that $g({\bf d})$ denote the number of graphs on
vertex set $[n]$ with degree sequence ${\bf d}$. We have
$$
\pr_{\mathcal{G}_{ {\bf d}}}(S,H)=\frac{g( S,[n]\setminus S, {\bf
d'})}{g({\bf d})}.
$$
The corollary now follows from the formulae for $g( S,[n] \setminus S, {\bf d'})$   in
Theorem~\ref{t:graphCount} and  $g({\bf d})$
 in Theorem~\ref{t:regularCount}.  \qed 

\ss

 Let $\pr_{\mathcal{G}_{n,d}}(S,H)$ denote the probability that $G_S=
H$ for a random $d$-regular graph $G$.
\begin{cor}\lab{c1:graphCount}
Given $0<s<n$, let $S=[s]\subset [n]$, let $H$ be a graph on
 vertex set $S$ with degree sequence
${\bf k}=(k_1,\ldots,k_s)$ with $k_i\le d$ for all $1\le i\le
s$, and put $h=\sum_{i=1}^h k_i$. Assume $ d=o((n-s)^{1/3})$. Then
\begin{eqnarray*}
\pr_{\mathcal{G}_{n,d}}(S,H)
&=&\exp\left(-\mu_0({\bf d}')-\mu_1({\bf d}')-\mu_2({\bf d}')+\frac{d^2-1}{4} +O(d^4/(dn-h))\right)\\
&&\hspace{.5cm}\times\prod_{i=1}^s[d]_{k_i}\frac{(dn-ds)!(dn/2)!2^{ds-h/2}}{((dn-2ds+h)/2)!(dn)!}\, ,
\end{eqnarray*}
where $d'_i=d-k_i$ for $i\in S$ and $d'_i=d$ for $i\notin S$, and  $\mu_i$ is defined as in Corollary~\ref{t2:graphCount}.
\end{cor}
\proof We apply    Corollary~\ref{t2:graphCount}. By the definition
of $\mu(\bf d)$, we immediately get that $\mu({\bf d})+\mu({\bf
d})^2=(d^2-1)/4$ when ${\bf d}$ is a constant sequence with each
term $d$. We also have $M({\bf d})=dn$, $M({\bf d}')=dn-h$,
$M_1({\bf d}',S)=ds-h$, $M_1({\bf d}',[n] \setminus S)=dn-ds$,   and
$d'_{\max}\le d$.  Moreover,
$$
\frac{ (d_{\max}')^4}{M({\bf
d}')}=\frac{d^4}{dn-h}=\frac{d^3}{n-h/d}\le\frac{d^3}{n-s}=o(1),
$$
since $h\le ds$ and $d=o((n-s)^{1/3})$. \qed \ss

The formula in Corollary~\ref{c1:graphCount} easily simplifies if the graph  $H$ is not too large.
\begin{cor}\lab{c2:graphCount}
Let $S$, $H$, ${\bf k}$ and $h$ be defined as in
Corollary~\ref{c1:graphCount}. If $ d=o(n^{1/3})$, $s^2d=o(n)$ and
$d^2s=o(n)$, then
\begin{eqnarray*}
&&\pr_{\mathcal{G}_{n,d}}(S,H)=\big(1+O((d^{3}+s^2d+d^2s)/n)\big)
(dn)^{-h/2}\prod_{i=1}^s [d]_{k_i}.
\end{eqnarray*}
\end{cor}

\proof Since
$d^2s=o(n)$, we have $h=O(ds)=o(n)$ and hence $ d^4/(dn-h)=O ( d^3/n).$
Similarly,
$$
M_1({\bf d}',R) = dn+O(ds),\quad M_i({\bf d}',L)=O(d^is)\ (i=1,2),\quad M_2({\bf d}',R)= d(d-1)(n-O(s))
$$
and hence from~\eqn{eq:mu0}--\eqn{eq:mu2},
$$
\mu_0({\bf d'})=\frac{d-1}{2}+O(ds/n), \quad
\mu_1({\bf d'}) =O\left(d^2s/n\right), \quad
\mu_2({\bf d'})=\frac{(d-1)^2}{4}+O(d^2s/n).
$$
Thus
$\mu_0({\bf d'})+ \mu_1({\bf d'})+\mu_2({\bf
d'})=(d^2-1)/4 +O(d^2s/n)$.
\ss

The corollary now follows upon applying Stirling's formula in the form $n!=\sqrt{2\pi n}(n/e)^n(1+O(n^{-1}))$ to obtain (ignoring negligible error terms)
$$
\frac{(dn-ds)!(dn/2)!2^{ds-h/2}}{((dn-2ds+h)/2)!(dn)!}=\left(\frac{dn}{e}\right)^{-h/2} \frac{(1-s/n)^{dn-ds} }{(1-2s/n+h/dn)^{(dn-2ds+h)/2}}. \qed
$$

Another interesting special case is when $H$ is empty.
\begin{cor}\lab{c3:graphCount}
Assume $ d=o(n^{1/3})$. Then for any $S\subset [n]$ with $
s=|S| <n/2$,
\begin{eqnarray*}
 \pr(S\ \mbox{is independent} )=  
\big(1+O(d^{3}/n)\big)\exp\left(f(d,\delta)\right)\prod_{i=1}^s\frac{(dn-ds)!(dn/2)!2^{ds}}{((dn-2ds)/2)!(dn)!},
\end{eqnarray*}
where $\delta=\delta(n)=s/n$, and
$$
f(d,\delta)= -\frac{\delta(d - 1) (  \delta d- 2 +\delta)}
 {4    (1-\delta)^2}.
 $$
\end{cor}
\proof This is a simple application of Corollary~\ref{c1:graphCount} with $h=0$, noting that
$$
\mu_0=\frac{(d-1)(n-2s)}{2(n-s)},\ \mu_1=\frac{(d-1)^2s}{2(n-s)},\ \mu_2=\frac{(d-1)^2(n-2s)^2}{4(n-s)^2}.\qed
$$
Note that if $d(n-2s)\to\infty$, then the probability that $S$ is independent under the conditions in Corollary~\ref{c3:graphCount} can be further simplified using Stirling's formula to
$$
 \big(1+O(d^{3}/n)+O(1/(dn-2ds))\big)
\sqrt{\frac{1-\delta}{1-2\delta}}\left(\frac{(1-\delta)^{1-\delta}}{(1-2\delta)^{(1-2\delta)/2}}\right)^{dn}\exp\left(f(d,\delta)\right).
$$

\section{The main switchings}
\lab{s:basics}


We can use the pairing model to generate B-graphs with the vertex
partition $L\cup R$ and the degree sequence ${\bf d}=\{d_1,\ldots,
d_n\}$.  Consider $n$ buckets representing the $n$ vertices. Let
each bucket $i$ contain $d_i$ points. Take a random pairing of
these points. 
We say a pairing is {\em restricted} if no pair has both ends in the
buckets representing vertices in $L$. Let $\mathcal{M}(L,R,{\bf d})$
be the class of all restricted pairings. Every such pairing
corresponds to a B-multigraph by contracting all points in each
bucket to form a vertex.  In the rest of the paper,  a bucket in a
pairing is also called a vertex.  A pair in a pairing is called a
{\em mixed (pure) pair} if it corresponds to a mixed (pure) edge in
the corresponding B-multigraph.  Thus, in a restricted pairing, each
pair is either mixed or pure; pure pairs have both points in a
vertex in $R$. Note that any simple B-graph corresponds to
$\prod_{i=1}^nd_i!$ restricted pairings in $\mathcal{M}(L,R,{\bf
d})$.
Hence, all simple B-graphs occur with the same probability in the
pairing model.


\ss

The main goal of this section is to compute the probability that a
B-multigraph generated by the pairing model is simple. We say that
$\{\{u_1,u'_1\},\{u_2,u'_2\},\{u_3,u'_3\}\}$ is a triple pair  if
$u_1$, $u_2$, $u_3$ are in one vertex and $u'_1$, $u'_2$, $u'_3$ are
in another vertex. We call the two vertices involved the {\em end
vertices} of the triple pair. If the end vertices are in $L$ and $R$
respectively,   the triple pair is called a {\em mixed triple pair},
and otherwise it is {\em pure}. Given a random restricted pairing,
let $T_1$ and $T_2$ be the number of mixed and pure triple
pairs respectively. 
In this section, there is only one degree sequence $\bf d$
referred to, so we drop the notation ${\bf d}$ from $M(\bf d)$ and
$M_i({\bf d},L)$, $M_i({\bf d},R)$, $\mu_i({\bf d})$ for simplicity.
Since $M_1(R)\ge M_1(L)$ by   assumption~\eqn{eq:assumption}, we
have $M_1(R)\ge M/2$.
\begin{lemma}\lab{l:triple} $\ex(T_1)=O(\d^4/M)$ and $\ex(T_2)=O(\d^4/M)$.
\end{lemma}
\proof For any two vertices $i\in L$ and $j\in R$, we compute the
probability that there is a triple pair with end vertices $i$ and
$j$.
 There are $\binom{d_i}{3}$ ways to
choose three points from the vertex $i$ and $\binom{d_j}{3}$ ways to
choose three points from the vertex $j$. There are $6$ ways to match
the six chosen points to form a triple pair. For
any positive even integer $m$, let $\mathscrC(m)$   denote the number of pairings of $m$
points. Then
$$
\mathscrC(m)=\prod_{i=0}^{m/2-1}(m-2i-1)=\frac{m!}{2^{m/2}(m/2)!}.
$$
 The probability for the three particular
pairs to occur is
$$
\frac{[M_1(R)-3]_{M_1(L)-3}\mathscrC(M_1(R)-M_1(L))}{[M_1(R)]_{M_1(L)}\mathscrC(M_1(R)-M_1(L))}\sim
M_1(R)^{-3}
$$
(noting that $M_1(R)\ge M_1(L)$ implies $M_1(R)\to\infty$).
 This is because the number
of ways to match the remaining $M_1(R)-3$ points in $L$ to points in
$R$, except for the three chosen points in the vertex $j$, is
$[M_1(R)-3]_{M_1(L)-3}$, and the number of matchings of the
remaining $M_1(R)-M_1(L)$ points in $R$ is
$\mathscrC(M_1(R)-M_1(L))$, whilst the total number of restricted
pairings is $[M_1(R)]_{M_1(L)}\mathscrC(M_1(R)-M_1(L))$. Hence we
have
\begin{eqnarray*}
\ex(T_1)&\sim&\sum_{i\in L}\sum_{j\in
R}6\binom{d_i}{3}\binom{d_j}{3}M_1(R)^{-3}=O\left(\left(\sum_{i\in
L}d_i^3\right)\left(\sum_{j\in
R}d_j^3\right)\right)M^{-3}\\
&=&O\left(\frac{\d^4 M_1(L) M_1(R)}{M^3}\right)=O\left(\frac{\d^4
}{M} \right),
\end{eqnarray*}
where the second equality uses $M/2\le M_1(R)\le M$.

A similar argument gives
\begin{eqnarray*}
\ex(T_2)&\sim&\sum_{i\in R}\sum_{j\in
R}6\binom{d_i}{3}\binom{d_j}{3}M_1(R)^{-3}=O\left(\left(\sum_{i\in R}d_i^3\right)\left(\sum_{j\in R}d_j^3\right)\right)M^{-3}\\
&=&O\left(\frac{\d^4 M_1(R)^2}{M^3}\right)=O\left(\frac{\d^4 }{M}
\right).\qed
\end{eqnarray*}

A pair
$\{u,u'\}$   is called a {\em loop} if $u$ and $u'$ are
contained in the same vertex and two pairs $\{u_1,u'_1\},
\{u_2,u'_2\}$ are called a {\em double pair}  if $u_1$, $u_2$ are in one
vertex and $u'_1$, $u'_2$ are in another vertex.
We call two loops that contain points from a common vertex  a {\em double
loop}. Let $I$ be the number of double loops. The proof of the following is a simple modification of the proof of the previous lemma, so is omitted.
\begin{lemma}\lab{l:doubleLoop}
$\ex(I)=O(\d^3/M)$. \qed
\end{lemma}

Lemmas~\ref{l:triple} and~\ref{l:doubleLoop} show that a.a.s.\ there
are no triple pairs or double loops in a random restricted pairing,
under the assumption $\d^4=o(M({\bf d}))$. So we only need to
consider   loops and double pairs.  In a restricted
pairing, there are two types of double pairs. One is that $u_1$,
$u_2$ are contained in a vertex in $L$ and $u'_1$, $u'_2$ are
contained in a vertex in $R$. The other is that all of $u_1$, $u_2$,
$u'_1$ and $u'_2$ are contained in vertices in $R$. We call the
former type {\em mixed} and the latter type {\em pure}.

 Let
$B_0$, $B_1$ and $B_2$ be the numbers of loops, mixed double pairs
and pure double pairs respectively. We first compute the expected
value of $B_i$ for $i=0,1,2$.  Recall
from~\eqn{eq:mu0}--\eqn{eq:mu2} that
\begin{eqnarray*}
\mu_0&=&\frac{(M_1(R)-M_1(L))M_2(R)}{2M_1(R)^2}, \ \
\mu_1=\frac{M_2(R)M_2(L)}{2M_1(R)^2},\ \ \mu_2=\mu_0^2.
\end{eqnarray*}
\begin{lemma}\lab{l:expectation}
For    $i=0,1,2$ we have
$ \ex B_i =O(\mu_i)$.  If
$\d=o(M^{1/3})$ and $M_1(R)-M_1(L)\to \infty$, then, more precisely, $\ex B_i\sim \mu_i$ for $i=0$ and $1$, and $\ex B_2=(1+o(1))\mu_2+o(1)$.
  \end{lemma}
\proof  Using small modifications of the proof  of
Lemma~\ref{l:triple},
we   immediately get
\begin{eqnarray*}
\ex B_0&=&\sum_{i\in
R}\binom{d_i}{2}\frac{[M_1(R)-2]_{M_1(L)}\mathscrC(M_1(R)-M_1(L)-2)}{[M_1(R)]_{M_1(L)}\mathscrC(M_1(R)-M_1(L))}\nonumber\\
&=&\sum_{i\in
R}\frac{[d_i]_2}{2}\frac{O\left(M_1(R)-M_1(L)\right)}{M_1(R)^2}=O(\mu_0); \nonumber\\
\end{eqnarray*}
\begin{eqnarray*}
\ex B_1&=&\sum_{i\in L}\sum_{j\in
R}2\binom{d_i}{2}\binom{d_j}{2}\frac{[M_1(R)-2]_{M_1(L)-2}\mathscrC(M_1(R)-M_1(L))}{[M_1(R)]_{M_1(L)}\mathscrC(M_1(R)-M_1(L))}\nonumber\\
&\sim&\frac{M_2(L)M_2(R)}{2}M_1(R)^{-2}=\mu_1;\nonumber\\
\end{eqnarray*}
\begin{eqnarray}
\ex B_2&=&\sum_{i,j\in R, i<j} 2\binom{d_i}{2}\binom{d_j}{2}\frac{[M_1(R)-4]_{M_1(L)}\mathscrC(M_1(R)-M_1(L)-4)}{[M_1(R)]_{M_1(L)}\mathscrC(M_1(R)-M_1(L))}\nonumber\\
&=&\frac{1}{2}\sum_{i\in R}\sum_{j\in
R}2\binom{d_i}{2}\binom{d_j}{2}\frac{[M_1(R)-4]_{M_1(L)}\mathscrC(M_1(R)-M_1(L)-4)}{[M_1(R)]_{M_1(L)}\mathscrC(M_1(R)-M_1(L))}\nonumber\\
&&-\frac{1}{2}\sum_{i\in R}2\binom{d_i}{2}\binom{d_i}{2}\frac{[M_1(R)-4]_{M_1(L)}\mathscrC(M_1(R)-M_1(L)-4)}{[M_1(R)]_{M_1(L)}\mathscrC(M_1(R)-M_1(L))}\lab{eq:doubleCount}\\
&=&\frac{M_2(R)^2}{4}\frac{O((M_1(R)-M_1(L))^2)}{M_1(R)^{4}}-\alpha
=O(\mu_2)-\alpha, \nonumber
\end{eqnarray}
where $\alpha =O(\d^3/M)$ is nonnegative. This gives the first part of the lemma.

If furthermore $\d=o(M^{1/3})$ and $M_1(R)-M_1(L)\to \infty$, then
all the $O(.)$ terms in the displayed equations above can be
replaced by $(1+o(1))(.)$. The lemma follows. \qed \ss

\begin{cor} \lab{c2:expectation} If $\d^4=o(M)$ and
$M_2(R)=O(\d^3)$, then the probability that there exists a loop or a
double pair is $O(\d^4/M)$.
\end{cor}
\proof If $\d^4=o(M)$ and $M_2(R)=O(\d^3)$, then $ \ex
B_0=O(M_2(R)/M_1(R))=O(\d^3 /M)$; $\ex
B_1=O(M_2(L)\d^3/M^2)=O(\d^4/M)$ (since $M_2(L)/M_1(R)\le
M_2(L)/M_1(L)\le \d$); $\ex B_2=O(\d^6 /M^2)=o(\d^2/M)$. The result
follows by the first moment principle.\qed \ss

We will need to prescribe some upper bounds on the likely values of
the random variables of interest. Define
$$
\eta(L)=M_2(L)/M_1(L), \quad \eta(R)=M_2(R)/M_1(R)
$$
and let
\bel{kdefs}
k_0=\max\{\ln M, 8\eta(L), 8\eta(R)\},\quad
k_1=k_2=\max\{\ln M, 8\eta(L)^2, 8\eta(R)^2\}\ (i=1,2).
\ee
Clearly
$\eta(L)=O(\d)$ and $\eta(R)=O(\d)$.
\begin{lemma} \lab{c1:expectation}
If $\d^4=o(M)$, then
$\pr\big(B_i\ge k_i\big)=O(M^{-1})$ for $i=0,1,2$.
\end{lemma}
\proof For any $h=o(\sqrt{M})$, the probability that there exist $h$
loops is bounded above by the $h$-th factorial moment of $B_0$.
Following the same pattern of proof  as for Lemma~\ref{l:triple},
this is at most
\begin{eqnarray}
&&\sum_{\substack{i_1,\ldots,i_h\in R\\i_1<\cdots<i_{h}}}\left(\prod_{j=1}^h\binom{d_{i_j}}{2}\right)\frac{[M_1(R)-2h]_{M_1(L)}\mathscrC(M_1(R)-M_1(L)-2h)}{[M_1(R)]_{M_1(L)}\mathscrC(M_1(R)-M_1(L))}\nonumber\\
&&\hspace{.2cm}\le \frac{M_2(R)^h}{2^hh!}\frac{[M_1(R)-2h]_{M_1(L)}\mathscrC(M_1(R)-M_1(L)-2h)}{[M_1(R)]_{M_1(L)}\mathscrC(M_1(R)-M_1(L))}\nonumber\\
&&\hspace{.2cm}=
\frac{M_2(R)^h}{2^hh!}\frac{\prod_{i=0}^{2h-1}(M_1(R)-M_1(L)-i)}{\prod_{i=0}^{2h-1}(M_1(R)-i)}
\left(\prod_{i=0}^{h-1}(M_1(R)-M_1(L)-2i-1)\right)^{-1}\nonumber\\
&&\hspace{.2cm}=\frac{M_2(R)^h}{2^hh!}\frac{\prod_{i=0}^{h-1}(M_1(R)-M_1(L)-2i)}{\prod_{i=0}^{2h-1}(M_1(R)-i)}\sim \frac{M_2(R)^h}{2^hh!}\frac{(M_1(R)-M_1(L))^h}{M_1(R)^{2h}}.\lab{tloop}
\end{eqnarray}
Since $M_1(R)=\Theta(M)$ and $h=o(\sqrt{M})$, this probability is at
most
\begin{eqnarray*}
\frac{M_2(R)^h}{2^hh!}\left(M_1(R)^h(1+o(1))\right)^{-1}\le
\left(\frac{eM_2(R)}{2hM_1(R)}\right)^h=\left(\frac{e\eta(R)}{2h}\right)^h.
\end{eqnarray*}
Similarly we have that for any $h=o(\sqrt{M})$, the probability that
there exist $h$ mixed double pairs is at most
\begin{eqnarray}
&&\sum_{\substack{i_1,\ldots,i_h\in L, j_1,\ldots,j_h\in R\\i_1<\cdots<i_{h}}}\left(\prod_{\ell=1}^h 2\binom{d_{i_{\ell}}}{2}\binom{d_{j_{\ell}}}{2}\right)\frac{[M_1(R)-2h]_{M_1(L)-2h}\mathscrC(M_1(R)-M_1(L))}{[M_1(R)]_{M_1(L)}\mathscrC(M_1(R)-M_1(L))}\nonumber\\
&&\hspace{.2cm}\le\frac{M_2(L)^hM_2(R)^h}{2^hh!}M_1(R)^{-2h}\le
\left(\frac{e}{2h}\cdot\frac{M_2(L)}{M_1(R)}\cdot\frac{M_2(R)}{M_1(R)}\right)^h,\lab{tdouble1}
\end{eqnarray}
and the probability that there exist $h$ pure double pairs is at
most
\begin{eqnarray}
&&\sum_{\substack{i_1,\ldots,i_h\in R, j_1,\ldots,j_h\in R\\i_1<\cdots<i_{h}}}\left(\prod_{\ell=1}^h 2\binom{d_{i_{\ell}}}{2}\binom{d_{j_{\ell}}}{2}\right)\frac{[M_1(R)-4h]_{M_1(L)}\mathscrC(M_1(R)-M_1(L)-4h)}{[M_1(R)]_{M_1(L)}\mathscrC(M_1(R)-M_1(L))}\nonumber\\
&&\hspace{.2cm}\le\frac{M_2(R)^hM_2(R)^h}{2^hh!}\frac{(M_1(R)-M_1(L))^{2h}}{M_1(R)^{4h}}.\lab{tdouble2}
\end{eqnarray}

Note that
$
\eta(L)$ and $\eta(R)$ are both bounded above by $\d$.
By the definition of $k_i$ in~\eqn{kdefs}, $k_i=O(\ln M+\d^2)$ for
  $i=0,1,2$. Since $\d^4=o(M)$, we therefore have $k_i=o(\sqrt{M})$. Hence
\begin{eqnarray*}
 \pr\big(B_0\ge k_0\big)&\le& \left(\frac{e\eta(R)}{2k_0}\right)^{k_0}\le
\left(\frac{e}{16}\right)^{\ln M}<M^{-1},\\
\pr\big(B_1\ge
k_1\big)&\le&\left(\frac{e}{2k_1}\cdot\eta(L)\cdot\eta(R)\right)^{k_1}\le\left(\frac{e}{16}\right)^{\ln
M}<M^{-1},\\
\pr\big(B_2\ge
k_2\big)&\le&\left(\frac{eM_2(R)^2}{2tM_1(R)^2}\right)^{k_2}=\left(\frac{e\eta(R)^2}{2k_2}\right)^{k_2}\le\left(\frac{e}{16}\right)^{\ln
M}<M^{-1}.\qed
\end{eqnarray*}
\begin{lemma}\lab{c:expectation} Assuming $\d^4=o(M)$,
\begin{description}
\item{(i)} if $M_2(R)=O(\d^5+\d^3\ln^2 M)$, then with probability $1-O(\d^4/M)$, $B_0\le\d+2$ and $B_i\le \d^2+2$ for all $i=1$ and $2$;
\item{(ii)} if $M_1(R)-M_1(L)=O(\d^4+\d^2\ln^2 M)$,
then with probability $1-O(\d^4/M)$, $B_0\le \d+2$ and $B_2\le \d^2+2$;
\item{(iii)} if $M_2(L)=O(\d^5+\d^3\ln^2 M)$, the with probability
$1-O(\d^4/M)$, $B_1\le \d^2+2$.
\end{description}
\end{lemma}
\proof These statements follow easily, after some simple
estimations, from~\eqn{tloop},~\eqn{tdouble1}
and~\eqn{tdouble2}.\qed\ss

We now redefine the values $k_i$ as follows.  Let $\zeta_0$,
$\zeta_1$ and $\zeta_2$ be (large) constants specified later. If
$M_2(R)\le \zeta_0(\d^5+\d^3\ln^2 M)$, use $k_0=\d+2$ and
$k_i=\d^2+2$ for $i=1$ and $2$; if $M_1(R)-M_1(L)\le
\zeta_1(\d^4+\d^2\ln^2 M)$, $k_0=\d+2$, and $k_2=\d^2+2$;
 if $M_2(L)\le \zeta_2(\d^5+\d^3\ln^2 M)$, use $k_1=\d^2+2$.
With the modified values, we have the following immediately from
the previous two results.
\begin{cor} \lab{c:mod}
 If $\d^4=o(M)$, then $\pr\big(B_i\ge
k_i\big)=O(\d^4/M)$ for $i=0,1,2$.
\end{cor}

 Define
$\mathC_{l_0,l_1,l_2}$ be the class of restricted pairings in
$\mathcal{M}(L,R,\bf d)$ that contains $l_0$ loops, $l_1$ mixed
double pairs, $l_2$ pure double pairs and no double loop or triple
pairs. Also, let $\pr( {\bf d})$ be the probability that a random
pairing $\PA\in \mathcal{M}(L,R,{\bf d})$ corresponds to a simple
B-graph.

The
 following corollary
is obtained from  Lemmas~\ref{l:triple} and~\ref{l:doubleLoop} and
 Corollary~\ref{c:mod} by noting that the sum of
$|\mathcal{C}_{l_0,l_1,l_2}|$ over {\em all} $l_0,l_1,l_2$
 is the total number of pairings with
$T_1=T_2=I=0$.
\begin{cor}\lab{c:probability}
\begin{eqnarray*}
\frac{1}{\pr({\bf
d})}=\big(1+O(\d^4/M)\big)\sum_{l_0=0}^{k_0}\sum_{l_1=0}^{k_1}\sum_{l_2=0}^{k_2}\frac{|\mathcal{C}_{l_0,l_1,l_2}|}{|\mathcal{C}_{0,0,0}|}.
\end{eqnarray*}
\end{cor}
With this corollary in mind, in the rest of the paper when
considering $|\mathcal{C}_{l_0,l_1,l_2}|$ we implicitly assume that
$0\le l_i\le k_i$ for  $i=0, 1$ and $2$.

Given a restricted pairing $\PA$,  we say the ordered pair of pairs
$((u_1,u'_1),(u_2,u'_2))$ forms a directed 2-path in $\PA$ if $u'_1$
and $u_2$ lie in the same vertex and the three vertices where $u_1$,
$u'_1$ and $u'_2$ lie in respectively are all distinct. We then say
that the two pairs $(u_1,u'_1)$ and $(u_2,u'_2)$ are adjacent. For
instance, the ordered pair of pairs $((1,2),(3,4))$ forms a directed
2-path in the four examples in Figure~\ref{f:2-paths}. Note that a
directed 2-path in a pairing corresponds to a directed 2-path in the
corresponding B-multigraph. Let $v$ denote the vertex where $u'_1$
and $u_2$ lie in. We say the directed 2-path
$((u_1,u'_1),(u_2,u'_2))$ in $\PA$ is {\em simple} if neither of
$\{u_1,u'_1\}$ and $\{u_2,u'_2\}$ is contained in a double pair and
there is no loop at $v$.

There are four types of directed 2-paths in which we are interested
in this paper. These 2-paths will be used later to define our
switching operations. Those with all vertices lying in $R$ are of
{\em type 1}.  A directed 2-path $((a,b),(c,d))$ is of {\em type 2}
if   $a$ lies in a vertex in $L$ and the other points all lie in
vertices in $R$, {\em type 3} if $a$ and $d$ are in vertices in  $L$
and the vertex containing $b$ and $c$ is in $R$, and {\em type 4} if
$a$ and $d$  lie in vertices in $R$ and the vertex containing $b$
and $c$ is in $L$.

  Given a restricted pairing $\PA$, let $t$
 be the
number of pure pairs in $\PA$. Then
\bel{tdef}
t=(M_1(R)-M_1(L))/2.
\ee
 Let $A_i(\PA)$ denote the number of
simple directed $2$-paths of type $i$ for $i=1,2,3,4$ and let
\bel{adef}
a_i(l_0,l_1,l_2)=\ex(A_i(\PA)\mid \PA\in\mathC_{l_0,l_1,l_2}).
\ee
Clearly $A_4(\PA)=\sum_{i\in
L}d(i)(d(i)-1)-O(l_1\d)=M_2(L)-O(l_1\d)$  for any
$\PA\in\mathC_{l_0,l_1,l_2}$ since the number of non-simple directed
2-path of type 4 is bounded by $O(l_1\d)$. \ss

  \begin{figure}[htb]
\vbox{\vskip .8cm
 \hbox{\centerline{\includegraphics[width=8cm]{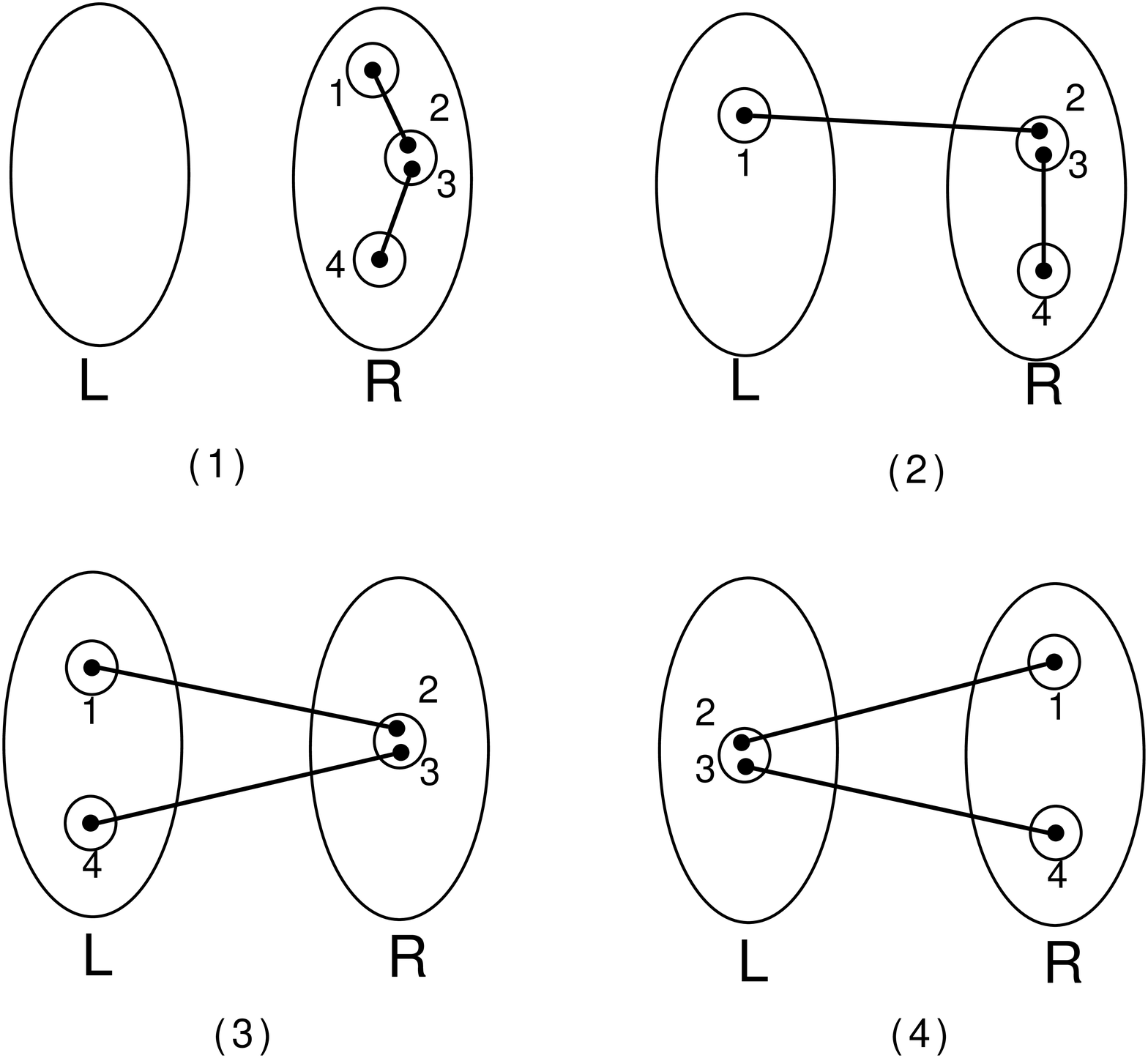}}}

\vskip .5cm \smallskip} \caption{\it  four different types of
$2$-paths}

\lab{f:2-paths}

\end{figure}

 The switching operations we are going to use are ideologically similar to
the switching operations used by McKay and Wormald~\cite{MW2}.
Although those switchings cannot be applied here because they do not
preserve the property of the pairings being restricted, they can easily be adjusted and adapted to our
current needs. The main twist is that there are a number of
alternative switchings available use, and we need to specify which
ones should be used, and for what values of the parameters, to
achieve the desired result. The following two switching operations
are used to prove Lemma~\ref{l:loop}.

\begin{description}
 \item{(a)} {\em $L_1$-switching}:  take a loop $\{2,3\}$ and two pure pairs $\{1,5\},\
 \{4,6\}$ such that the six
 points are located in the five distinct vertices  as drawn in
 Figure~\ref{f:swt1}.
 Replace the three pairs $\{2,3\},\{1,5\},\{4,6\}$ by
 $\{1,2\},\{3,4\},\{5,6\}$.

\item{(b)} {\em $L_2$-switching}:  take a loop $\{2,3\}$ and two mixed pairs
$\{1,5\},\
 \{4,6\}$  such that the six
 points are located in the five distinct vertices  as drawn in
 Figure~\ref{f:at-swt1}.
 Replace the three pairs $\{2,3\},\{1,5\},\{4,6\}$ by
 $\{1,2\},\{3,4\},\{5,6\}$.
\end{description}

\begin{figure}[htb]
\vbox{\vskip .8cm
 \hbox{\centerline{\includegraphics[width=6cm]{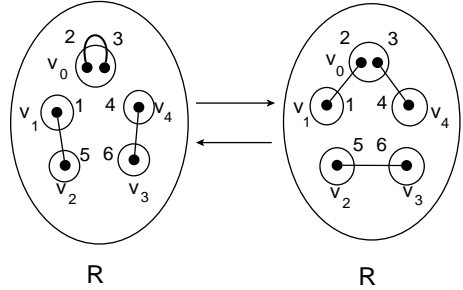}}}

\vskip .5cm \smallskip} \caption{\it  $L_1$-switching}

\lab{f:swt1}

\end{figure}

\begin{figure}[htb]
\vbox{\vskip .8cm
 \hbox{\centerline{\includegraphics[width=9cm]{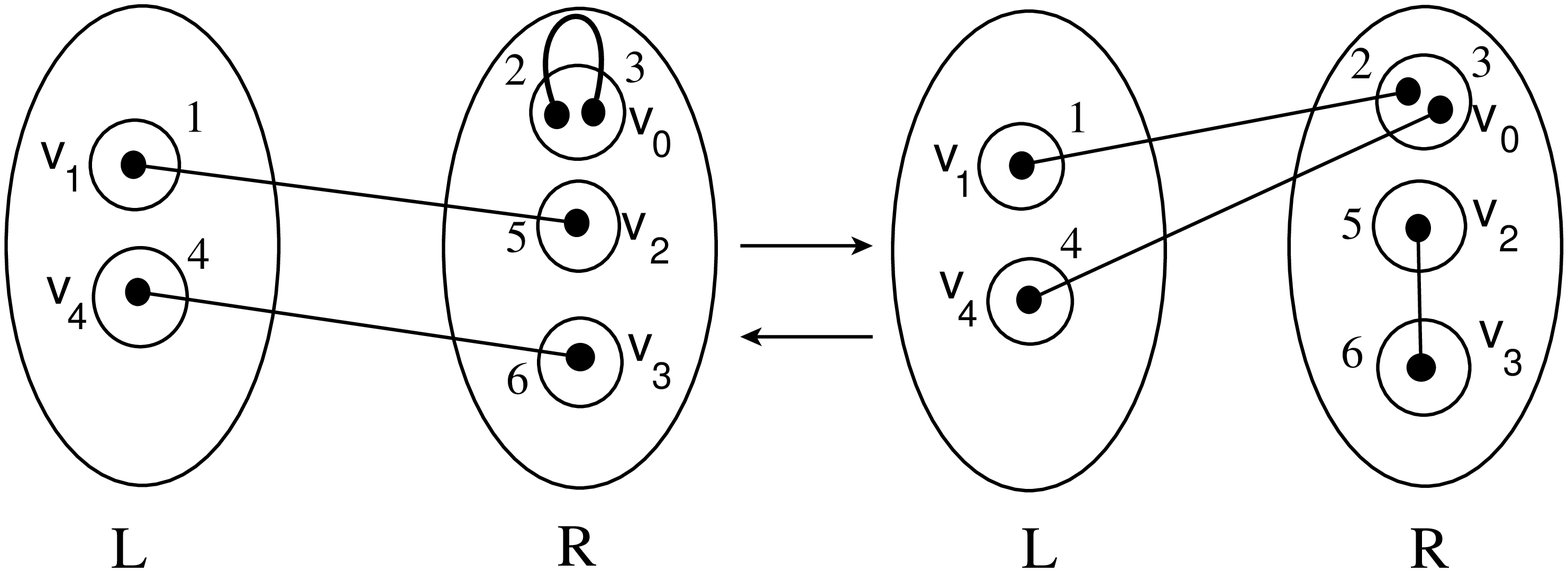}}}

\vskip .5cm \smallskip} \caption{\it  $L_2$-switching}

\lab{f:at-swt1}

\end{figure}

For any switching operation that converts a pairing $\PA_1$ to
another pairing $\PA_2$, we call the operation that converts $\PA_2$
to $\PA_1$ the inverse of that switching. Thus, the {\em inverse
$L_1$-switching} can be defined as follows. Take a 2-directed path
(not necessarily simple) $((1,2),(3,4))$ of type 1 and a pure pair
$\{5,6\}$ such that the points $1$, $2$, $4$, $5$ and $6$ lie in
five distinct vertices. Replace $\{1,2\}$, $\{3,4\}$ and $\{5,6\}$
by $\{2,3\}$, $\{1,5\}$ and $\{4,6\}$. The inverse $L_2$-switching
can be defined in the same way.

The following lemma
  estimates the ratio
$|\mathC_{l_0,l_1,l_2}|/|\mathC_{l_0-1,l_1,l_2}|$  by counting ways
to perform certain $L_1$-switchings and their inverses. We express
the present results in terms of the numbers $a_i(l_0,l_1,l_2)$,
defined in~\eqn{adef}, whose estimation we postpone till later.
\begin{lemma}\lab{l:loop}
Let $a_1=a_1(l_0-1,l_1,l_2)$ and $a_3=a_3(l_0-1,l_1,l_2)$. Assume
$l_0\ge 1$. Then
\begin{eqnarray*}
(i):&& \mbox{If}\ t\ge 1,\\
&&\frac{|\mathC_{l_0,l_1,l_2}|}{|\mathC_{l_0-1,l_1,l_2}|}=\frac{a_1}{4l_0t}(1+O(\d
^2/t+(l_0+l_2)/t)), \\
 (ii):&& \mbox{If}\ M_1(L)\ge 1\ \mbox{and}\ t\ge 1,\\
&&\frac{|\mathC_{l_0,l_1,l_2}|}{|\mathC_{l_0-1,l_1,l_2}|}=\frac{ta_3}{l_0M_1(L)^2}(1+O(\d
^2/M_1(L)+\d ^2/t+l_1/M_1(L)+(l_0+l_2)/t)).
\end{eqnarray*}
\end{lemma}

\proof Let $\PA\in\mathC_{l_0,l_1,l_2}$ and we consider the number
of $L_1$-switching operations that convert $\PA$ to some $\PA'\in
\mathC_{l_0-1,l_1,l_2}$. For the purpose of counting, we label the
points in the pairs that are under consideration  as shown in
Figure~\ref{f:swt1}. So for any pair under consideration,  we will
incorporate in our counting the number of ways we can label the
points in the pair. Let $N$ denote the number of ways to choose the
pairs and label the points in them so that an $L_1$-switching can be
applied to these pairs, which converts $\PA$ to some $\PA'\in
\mathC_{l_0-1,l_1,l_2}$ without any simultaneously created loops or
double pairs. This implies that the switching operations counted by
$N$ destroy only one loop and there is no simultaneous creation or
destruction of other loops or double pairs.

 We first give a rough count of $N$, that includes some
forbidden cases (due to creating double pairs, etc) and then
estimate the error. There are $l_0$ ways to choose a loop $e_0$ and
$t(t-1)$ ways to choose $(e_1,e_2)$, an ordered pair of two distinct
pure pairs. For any chosen loop $e_0$, there are two ways to
distinguish the two end points to label the points $2$ and $3$ as
shown in Figure~\ref{f:swt1}. For each of the other pairs, there are
two ways to label its two endpoints, as $1$ and $5$, or $4$ and $6$,
as the case may be. Hence a rough estimation of $N$ is $8l_0t(t-1)$,
including the count of some forbidden choices of $e_0$, $e_1$ and
$e_2$, which we estimate next. Let the vertices that contain points
$2,1,5,6,4$ be denoted by $v_0,v_1,v_2,v_3,v_4$ respectively as
shown in Figure~\ref{f:swt1}. The only possible exclusions caused by
invalid choices in the above are the following:
\begin{description}
\item{(a)} the loop $e_0$ is adjacent to $e_1$ or $e_2$, or $e_1$ is
adjacent to $e_2$, in which case, the $L_1$-switching is not
applicable since the definition of the $L_1$-switching excludes
cases where the edges are adjacent because it requires the end
vertices to be distinct;

\item{(b)} there exists a pair between $\{v_0,v_1\}$, or $\{v_0,v_4\}$, or $\{v_2,v_3\}$  in
$\PA$, in which case there will be more double pairs created after
the $L_1$-switching is applied;

\item{(c)} the pair $e_1$ or $e_2$ is a loop or is contained in a
double pair, in which case there is a simultaneously destroyed loop
or double pair.

\end{description}

First we show that the number of exclusions from case (a) is
$O(l_0t\d)$. The number of pairs of $(e_0,e_1)$ is at most $l_0t$.
For any given $e_0$ and $e_1$, the number of ways to choose a pair
$e_2$ such that $e_2$ is adjacent to $e_0$ or $e_1$ is at most $2\d$
since both $e_0$ and $e_2$ are adjacent to at most $\d$ pairs. Hence
the number of triples of $(e_0,e_1,e_2)$ such that $e_2$ is adjacent
to either $e_0$ or $e_1$ is at most $2l_0t\d$. By symmetry,  the
number of triples of $(e_0,e_1,e_2)$ such that $e_1$ is adjacent to
either $e_0$ or $e_2$ is also at most $2l_0t\d$. Hence the number of
exclusions from case (a) is $O(l_0t\d)$.

Next we show that the number of exclusions from case (b) is
$O(l_0t\d ^2)$. As just explained, the number of pairs of
$(e_0,e_1)$ is at most $l_0t$. For any given $e_0$ and $e_1$, the
number of ways to choose a pair $e_2$ such that $v_3$ is adjacent to
$v_2$ or $v_4$ is adjacent to $v_0$ is at most $2\d^2$, since both
$e_0$ and $e_1$ have at most $\d^2$ edges that are of distance $2$
away. Hence the number of triples $(e_0,e_1,e_2)$ such that $v_3$ is
adjacent to $v_2$ or $v_4$ is adjacent to $v_0$ is $O(l_0t\d^2)$. By
symmetry, the number of triples $(e_0,e_1,e_2)$ such that $v_3$ is
adjacent to $v_2$ or $v_0$ is adjacent to $v_1$ is $O(l_0t\d^2)$.
Hence the number of exclusions from case (b) is $O(l_0t\d ^2)$.

Now we show that the number of exclusions from case (c) is
$O(l_0^2t+l_0tl_2)$. The number of ways to choose $e_0,e_1,e_2$ such
that $e_1$ or $e_2$ is a loop is at most $2l_0^2t$ and the number of
ways to choose these three pairs such that $e_1$ or $e_2$ is
contained in a double pair is at most $2\cdot l_0t\cdot
2l_2=O(l_0tl_2)$. Hence the number of exclusions from case (c) is
$O(l_0^2t+l_0tl_2)$.

Thus, the number of exclusions in the calculation of $N$ is
 $O(l_0t\d ^2+l_0^2t+l_0tl_2)$. So
$N=8l_0t^2(1+O(\d ^2/t+(l_0+l_2)/t))$.

Now choose an arbitrary pairing $\PA'\in \mathC_{l_0-1,l_1,l_2}$.
Let $N'$ be the number of ways to choose the pairs and label points
in them so that an inverse $L_1$-switching operation can be applied
to these pairs such that $\PA'$ is converted to some $\PA\in
\mathC_{l_0,t_1,t_2}$ without any simultaneously destroyed loops or
double pairs.  To apply this operation we need to choose
$e'_0,e'_1,e'_2$, such that $(e'_0,e'_1)$ is a simple directed
2-path of type 1 and $e'_2$ is a pure pair. We consider the directed
2-path $(e'_0,e'_1)$ because it automatically gives a unique way of
distinguishing vertices $v_1$, $v_0$ and $v_4$ and labelling points
as $1$, $2$, $3$ and $4$ in Figure~\ref{f:swt1}. There are
$A_1(\PA')$ simple directed $2$-paths of type 1, and hence
$A_1(\PA')$ ways to choose the points as $1$, $2$, $3$ and $4$. The
number of ways to choose a pure pair $e'_2$ is $t$ and so there are
$2t$ ways to fix the vertices $v_2$, $v_3$ and the points $\{5,6\}$.
The only possible exclusions to the above choices are listed the
following cases.
\begin{description}
\item{(a)} There
exists a pair between $\{v_1,v_2\}$ or $\{v_3,v_4\}$ in $\PA'$,
since then more double pairs will be created if the inverse
$L_1$-switching is applied.

\item{(b)} The pair  $e'_2$ is a loop, in which case the inverse $L_1$-switching is not applicable, or $e'_2$ is
contained in a double pair, in which case a double pair is destroyed
after the application of the inverse $L_1$-switching.


\item{(c)} The pair $e_2'$ is adjacent to the 2-path or is contained in the 2-path, in which case the inverse $L_1$-switching operation is not applicable.

\end{description}

The number of exclusions from case (a) is
   $O(A_1(\PA')\d^2)$ and the numbers of exclusions from case (b) and (c)
   are
   $O(A_1(\PA')l_0+A_1(\PA')l_2)$ and $O(A_1(\PA')\d)$
   respectively.

   Thus, the number of exclusions from case (a)--(d) is
   $O(A_1(\PA')\d^2+A_1(\PA')l_0+A_1(\PA')l_2)$.
 So
 $$
 \ex(N')=\ex\big(2A_1t(1+O(\d^2 /t+(l_0+l_2)/t))\mid \PA'\in\mathC_{l_0-1,l_1,l_2}\big)=2a_1t(1+O(\d^2 /t+(l_0+l_2)/t)).
 $$
We count the pairs of $(\PA,\PA')$ such that $\PA\in
\mathC_{l_0,t_1,t_2}$, $\PA'\in \mathC_{l_0-1,l_1,l_2}$, and $\PA'$
is obtained by applying an $L_1$-switching to $\PA$, which destroys
only one loop without any simultaneously created loops or double
pairs. Then the number of such pairs of pairings equals to both
$|\mathC_{l_0,l_1,l_2}|\ex(N)$ and
$|\mathC_{l_0-1,l_1,l_2}|\ex(N')$. Thus,

$$
\frac{|\mathC_{l_0,l_1,l_2}|}{|\mathC_{l_0-1,l_1,l_2}|}=\frac{a_1}{4l_0t}(1+O(\d
^2/t+(l_0+l_2)/t)).
$$
This proves part (i) of Lemma~\ref{l:loop}. Analogously we can
deduce the following by analysing the $L_2$-switching and its
inverse.
\begin{eqnarray*}
\frac{|\mathC_{l_0,l_1,l_2}|}{|\mathC_{l_0-1,l_1,l_2}|}&=&\frac{2ta_3+O(\d
^2 a_3)+O(l_0a_3+l_2a_3)}{2l_0M_1(L)^2+O(\d
^2M_1(L)l_0+l_0M_1(L)l_1)}\\
&=&\frac{ta_3}{l_0M_1(L)^2}(1+O(\d ^2/M_1(L)+\d
^2/t+(l_0+l_2)/t+l_1/M_1(L))).
\end{eqnarray*}
Then we obtain part (ii) of Lemma~\ref{l:loop}. \qed

\ss

We use the following two switching operations to prove
Lemma~\ref{l:double1}.

\begin{description}
\item{(a)} {\em $D_1$-switching}:  take a mixed double pair
$\{\{3,4\},\{5,6\}\}$ and also two pure pairs $\{1,2\}$ and
$\{7,8\}$ such that the eight points are located in the six distinct
vertices as shown in Figure~\ref{f:swt2}. Replace the four pairs by
$\{1,3\},\{5,7\},\{2,4\},\{6,8\}$.  \ss

\begin{figure}[htb]
\vbox{\vskip .8cm
 \hbox{\centerline{\includegraphics[width=9cm]{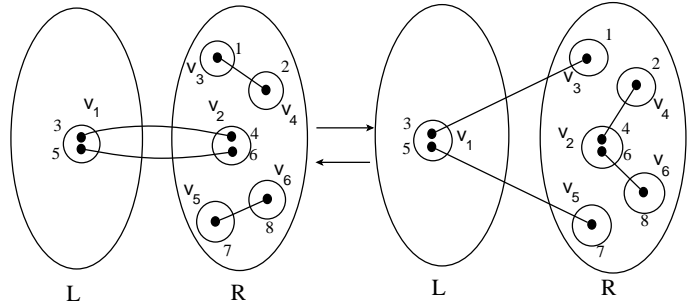}}}

\vskip .5cm \smallskip} \caption{\it $D_1$-switching}

\lab{f:swt2}

\end{figure}

\item{(b)} {\em $D_2$-switching}: take a mixed double pair $\{\{3,4\},\{5,6\}\}$
and also two mixed pairs $\{1,2\}$ and $\{7,8\}$ such that the eight
points are located in the six distinct vertices as shown in
Figure~\ref{f:at-swt2}. Replace the four pairs by
$\{1,4\},\{6,7\},\{2,3\},\{5,8\}$.
\end{description}

\begin{figure}[htb]
\vbox{\vskip .8cm
 \hbox{\centerline{\includegraphics[width=8cm]{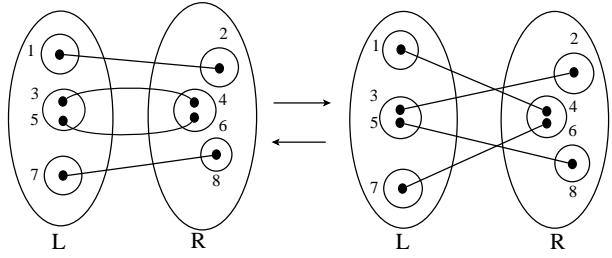}}}

\vskip .5cm \smallskip} \caption{\it  $D_2$-switching}

\lab{f:at-swt2}

\end{figure}

 The
inverse switchings are defined analogously to the earlier ones. For
instance, the inverse $D_1$-switching is defined as follows. Take a
directed $2$-path $((1,3),(5,7))$ of type $4$ and a directed
$2$-path $((2,4),(6,8))$ of type $1$ such that the eight points are
located in six distinct vertices as shown in Figure~\ref{f:swt2}.
Replace these four pairs by $\{1,2\}$, $\{3,4\}$, $\{5,6\}$  and
$\{7,8\}$.
\begin{lemma}\lab{l:double1}
Let $a_1=a_1(0,l_1-1,l_2)$ and $a_3=a_3(0,l_1-1,l_2)$. Assume
$l_1\ge 1$. Then
\begin{eqnarray*}
(i):&&\mbox{If}\ t\ge 1\ \mbox{and}\ M_2(L)\ge 1,\\
&&\frac{|\mathC_{0,l_1,l_2}|}{|\mathC_{0,l_1-1,l_2}|}=\frac{M_2(L)a_1}{8l_1t^2}(1+O(\d
^3/M_2(L)+\d ^2/t+l_2/t+l_1\d/M_2(L)));\\
(ii):&&\mbox{If}\ M_1(L)\ge 1\ \mbox{and}\ M_2(L)\ge 1,\\
&&\frac{|\mathC_{0,l_1,l_2}|}{|\mathC_{0,l_1-1,l_2}|}=\frac{a_3M_2(L)}{2l_1M_1(L)^2}(1+O(\d
^2/M_1(L)+\d ^3/M_2(L)+l_1/M_1(L)+l_1\d/M_2(L))).
\end{eqnarray*}

\end{lemma}

\proof For a given pairing $\PA\in\mathC_{0,l_1,l_2}$, let $N$ be
the number of ways to choose the pairs and label the points in them
so that a $D_1$-switching can be applied to these pairs such that
$\PA$ is converted to some $\PA'\in \mathC_{0,l_1-1,l_2}$ without
simultaneously creating any loops and double pairs. In order to
apply  a $D_1$-switching operation, we need to choose a mixed double
pair $\{e_1,e_2\}$ and an ordered pair of distinct pure pairs
$(e_3,e_4)$.
  The number of ways to choose
$\{e_1,e_2\}$ is $l_1$ in $\mathC_{0,l_1,l_2}$ and hence the number
of ways to label the points as $3,4,5,6$ is $2l_1$. The number of
ways to choose the ordered pair of pure pairs $(e_3,e_4)$ is
$t(t-1)$. For any chosen $(e_3,e_4)$, there are 4 ways to label
points as $1,2,7,8$. Let the vertices that contain points
$3,4,1,2,7,8$ be $v_1,v_2,v_3,v_4,v_5,v_6$ as shown in
Figure~\ref{f:swt2}. Hence a rough count of $N$ is $8l_1t(t-1)$
including the count of a few forbidden choices of $e_1,e_2,e_3,e_4$,
which are listed as follows.
\begin{description}
\item{(a)} The pair $e_1$ is adjacent to $e_3$ or $e_4$, or $e_3$ is
adjacent to $e_4$, in which case the $D_1$-switching is not
applicable.

\item{(b)} There exists a pair between $\{v_1,v_3\}$, or $\{v_2,v_4\}$, or $\{v_2,v_6\}$, or
$\{v_1,v_5\}$ in $\PA$, since another double pair will be created
after the $D_1$-switching is applied.

\item{(c)} The pair $e_3$ or $e_4$ is contained in a
double pair, since another double pair is destroyed after the
$D_1$-switching is applied.
\end{description}
The numbers of forbidden choices of $e_1,e_2,e_3,e_4$ coming from
case (a), (b) and (c) are $O(l_1t\d)$, $O(l_1t\d^2)$ and
$O(l_1tl_2)$ respectively.  So $N=8l_1t^2(1+O(\d^2 /t+l_2/t))$.

For a given pairing $\PA'\in\mathC_{0,l_1-1,l_2}$, let $N'$ be the
number of ways to choose the pairs and label the points in them so
that an inverse $D_1$-switching operation can be applied to these
pairs which converts $\PA'$ to some $\PA\in\mathC_{0,l_1,l_2}$
without destroying any loops or double pairs simultaneously. In
order to apply such an operation, we need to choose two simple
directed 2-paths, one of type 1 and the other of type 4. There are
$A_1(\PA')$ simple directed 2-paths of type 1, each of which gives a
way of labelling points as $2,4,6,8$, and there are $A_4(\PA')$
simple directed 2-paths of type 4, each of which gives a way of
labelling points as $1,3,5,7$. Hence a rough count of $N'$ is
$A_1(\PA')A_4(\PA')$ including the counts of a few forbidden choices
of such two 2-paths which are listed in the following two cases.
\begin{description}
\item{(a)} If we have $v_i=v_j$, for $i\in \{3,5\}$ and
$j\in\{2,4,6\}$, then the operation is not applicable.

\item{(b)} If there already exists a
pair between $\{v_1,v_2\}$, or $\{v_3,v_4\}$, or $\{v_5,v_6\}$ in
$\PA'$, then more than one double pair will be created in this case
if the inverse $D_1$-switching is applied.
\end{description}
The numbers of forbidden choices of the two directed 2-paths from
case (a) and (b) are respectively $O(A_1(\PA')\d^2)=O(a_1\d ^2)$ and
$O(A_1(\PA')\d^3)=O(a_1\d ^3)$. So
$\ex(N')=\ex\big(A_1(\PA')A_4(\PA')\mid
\PA'\in\mathC_{0,l_1-1,l_2}\big)+O(a_1\d
^3)=a_1(M_2(L)-O(l_1\d))(1+O(\d ^3/M_2(L)))$. Since $l_1\ge 1$, we
have $M_2(L)\ge 1$. Hence
\begin{eqnarray*}
\frac{|\mathC_{0,l_1,l_2}|}{|\mathC_{0,l_1-1,l_2}|}&=&\frac{a_1M_2(L)(1+O(\d
^3/M_2(L))+O(l_1\d/M_2(L)))}{8l_1t^2(1+O(\d ^2/t)+O(l_2/t))}\\
&=&\frac{a_1M_2(L)}{8l_1t^2}(1+O(\d ^3/M_2(L)+\d
^2/t+l_2/t+l_1\d/M_2(L))),
\end{eqnarray*}
and this shows part (i) of Lemma~\ref{l:double1}. Similarly we can
obtain part (ii) by analysing the $D_2$-switching and its inverse.
\qed

\ss

The following two switching operations are used for the next lemma.
\begin{description}
\item{(a)} {\em $D_3$-switching}:  take a pure double pair $\{\{1,2\},\{3,4\}\}$
and also two pure pairs $\{5,6\}$ and $\{7,8\}$ such that the eight
points are located in the six distinct vertices as shown in
Figure~\ref{f:swt3}. Replace the four pairs by
$\{1,5\},\{2,6\},\{3,7\},\{4,8\}$.  \ss

\begin{figure}[htb]
\vbox{\vskip .8cm
 \hbox{\centerline{\includegraphics[width=6cm]{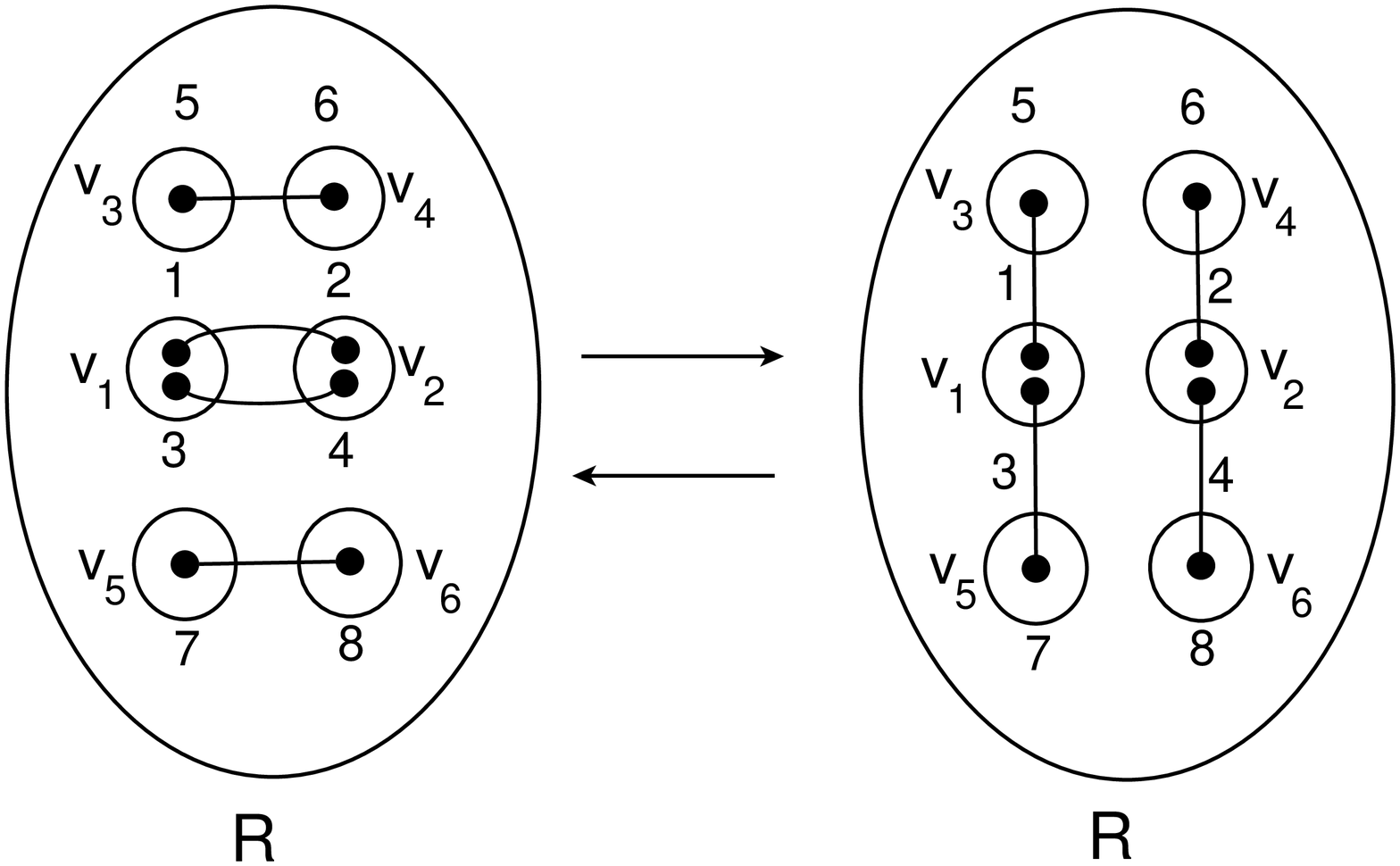}}}

\vskip .5cm \smallskip} \caption{\it  $D_3$-switching}

\lab{f:swt3}

\end{figure}

\smallskip

\item{(a)} {\em $D_4$-switching}:  take a pure double pair $\{\{1,2\},\{3,4\}\}$
and also four mixed pairs $\{5,6\}$, $\{7,8\}$, $\{9,10\}$,
$\{11,12\}$ such that the twelve points are located in the ten
distinct vertices as shown in Figure~\ref{f:at-swt3}. Replace the six
pairs by $\{6,10\},\{8,12\},\{1,5\},\{3,9\}$,$\{2,11\}$, $\{4,7\}$.
\end{description}

\begin{figure}[htb]
\vbox{\vskip .8cm
 \hbox{\centerline{\includegraphics[width=8cm]{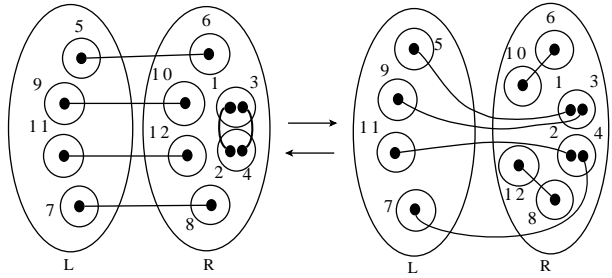}}}

\vskip .5cm \smallskip} \caption{\it  $D_4$-switching}

\lab{f:at-swt3}

\end{figure}

The inverse switchings are defined in the obvious way. For example,
for the inverse  of the $D_3$-switching, take two directed paths of
type 1, $((5,1),(3,7))$ and $((6,2),(4,8))$, such that the eight
points are located in six distinct vertices as shown in
Figure~\ref{f:swt3}. Replace these four pairs by $\{5,6\}$,
$\{1,2\}$, $\{3,4\}$, $\{7,8\}$.  Define $b_i(l_0,l_1,l_2)=\ex(A_i(\PA)^2\mid
\PA\in\mathC_{l_0,l_1,l_2})$ for $i=1$ and $3$.
\begin{lemma}\lab{l:double2}
 Assume $l_2\ge
1$.  For $i=1,3$, let $b_i=b_i(0,0,l_2-1)$ for short. Then
\begin{eqnarray*}
(i):&&\mbox{If}\ t\ge 1\ \mbox{and}\ b_1\ge 1,\\
&&\frac{|\mathC_{0,0,l_2}|}{|\mathC_{0,0,l_2-1}|}=\frac{b_1}{16l_2t^2}(1+O(\d
^2/t+\d ^3a_1/b_1+l_2/t)).\\
(ii):&&\mbox{If}\ M_1(L)\ge 1,\ b_3\ge 1 \ \mbox{and}\ t\ge 1,\\
&&\frac{|\mathC_{0,0,l_2}|}{|\mathC_{0,0,l_2-1}|}=\frac{t^2
b_3}{l_2M_1(L)^4}(1+O(\d ^3a_3/b_3+\d ^2/M_1(L)+l_2/t)).
\end{eqnarray*}

\end{lemma}

\proof For a given pairing $\PA\in\mathC_{0,0,l_2}$, let $N$ be the
number of ways to choose the pairs and label the points in them so
that a $D_3$-switching operation can be applied, which converts
$\PA$ to some $\PA'\in\mathC_{0,0,l_2-1}$ without creating any loops
and double pairs simultaneously. In order to apply a $D_3$-switching
operation, we need to choose a pure double pair $\{e_1,e_2\}$ and an
ordered pair of distinct pure pairs $(e_3,e_4)$. The number of ways
to choose $\{e_1,e_2\}$ is $l_2$ in $\mathC_{0,0,l_2}$ and there are
four ways to label the points as $1,2,3,4$ for any chosen double
pair. The number of ways to choose an ordered pair of two pure pairs
$(e_3,e_4)$ is $t(t-1)$ and hence the number of ways to label the
points as $5,6,7,8$ is $4t(t-1)$. Hence a rough count of $N$ is
$16l_2t(t-1)$ including the counts of forbidden choices of pairs
$e_1,\ldots,e_4$ which we estimate next. Let the vertices that
contain points $1,2,5,6,7,8$ be $v_1,v_2,v_3,v_4,v_5,v_6$ as shown
in Figure~\ref{f:swt3}. The forbidden choices of the pairs
$e_1,\ldots,e_4$ are listed in the following three cases.
\begin{description}
\item{(a)} When $e_1$ is adjacent to $e_3$ or $e_4$ or when $e_3$
is adjacent to $e_4$, then the $D_3$-switching is not applicable.

\item{(b)} If there exists a pair between $\{v_1,v_3\}$, or $\{v_2,v_4\}$, or $\{v_1,v_5\}$, or
$\{v_2,v_6\}$ in $\PA$, then more double pairs will be created after
the application of the switching operation.

\item{(c)} If $e_3$ or $e_4$ is contained in a
double pair, then another double pair would be destroyed after the
application of the switching operation.

\end{description}

 The numbers of
forbidden choices of $e_1,\ldots,e_4$ coming from (a),(b) and (c)
are $O(l_2t\d)$, $O(l_2t\d ^2)$ and $O(l_2^2t)$ respectively. So
$N=16l_2t^2(1+O(\d ^2/t+l_2/t))$.

For any pairing $\PA'\in\mathC_{0,0,l_2-1}$, let $N'$ be the number
of ways to choose the pairs and label the points in them so that an
inverse $D_3$-switching can be applied to these pairs, which
converts $\PA'$ to some $\PA\in\mathC_{0,0,l_2}$ without
simultaneously destroying any loops or double pairs. In order to
apply such an operation, we need to choose an ordered pair of
distinct simple directed 2-paths of type 1. The number of ways to do
that is $A_1(\PA')(A_1(\PA')-1)$. So the number of ways to label the
points $1,2,\ldots, 8$ is $A_1(\PA')(A_1(\PA')-1)$, which gives a
rough count of $N'$. The forbidden choices of the two paths whose
counts should be excluded from $N'$ are listed in the following
cases.
\begin{description}
\item{(a)} The two paths share some common vertex or common pair. In this case the
inverse $D_3$-switching is not applicable.
\item{(b)} There exists a pair between $\{v_1,v_2\}$ or $\{v_3,v_4\}$ or $\{v_5,v_6\}$  in $\PA'$. In this case, more double pairs will be created after
the inverse $D_3$-switching operation is applied.
\end{description}
The numbers of ways to choose the ordered pair of 2-paths in case
(a) and (b) are $O(A_1(\PA')\d^2)$ and $O(A_1(\PA')\d ^3)$
respectively. Thus, $\ex(N')=b_1(1+O(\d ^3a_1/b_1))$.

Hence
$$
\frac{|\mathC_{0,0,l_2}|}{|\mathC_{0,0,l_2-1}|}=\frac{b_1}{16l_2t^2}(1+O(\d
^2/t+\d ^3a_1/b_1+(l_0+l_2)/t)).
$$

Similarly by analysing the $D_4$-switching and its inverse, we
obtain Lemma~\ref{l:double2}(ii).\qed

\section{More switchings to estimate $a$'s and $b$'s}
\lab{s:moreswitchings} The lemmas in the previous section give
ratios of the sizes of `adjacent' classes $\mathC_{i,j,k}$, but
those estimates are in terms of   $a_i(l_0,l_1,l_2)$ ($i=1,2,3$)
defined in~\eqn{adef},   $b_i$ ($i=1,3$) defined just before
Lemma~\ref{l:double2}, and $t$ defined in~\eqn{tdef}. In this
section, we use further switchings to estimate the  values of these
variables.
  The following two switchings
are used for $a_i$.
\begin{description}
\item{(a)} {\em $S_1$-switching}: Take a mixed pair and label the points in it by
$\{1,2\}$ as shown in Figure~\ref{f:swt4}. Take a simple directed
$2$-path that is  vertex disjoint from the chosen mixed pair. Label
the points by $3,4,5,6$. Replace these three pairs by $\{2,3\}$,
$\{1,4\}$ and $\{5,6\}$.

\begin{figure}[htb]
\vbox{\vskip .8cm
 \hbox{\centerline{\includegraphics[width=9cm]{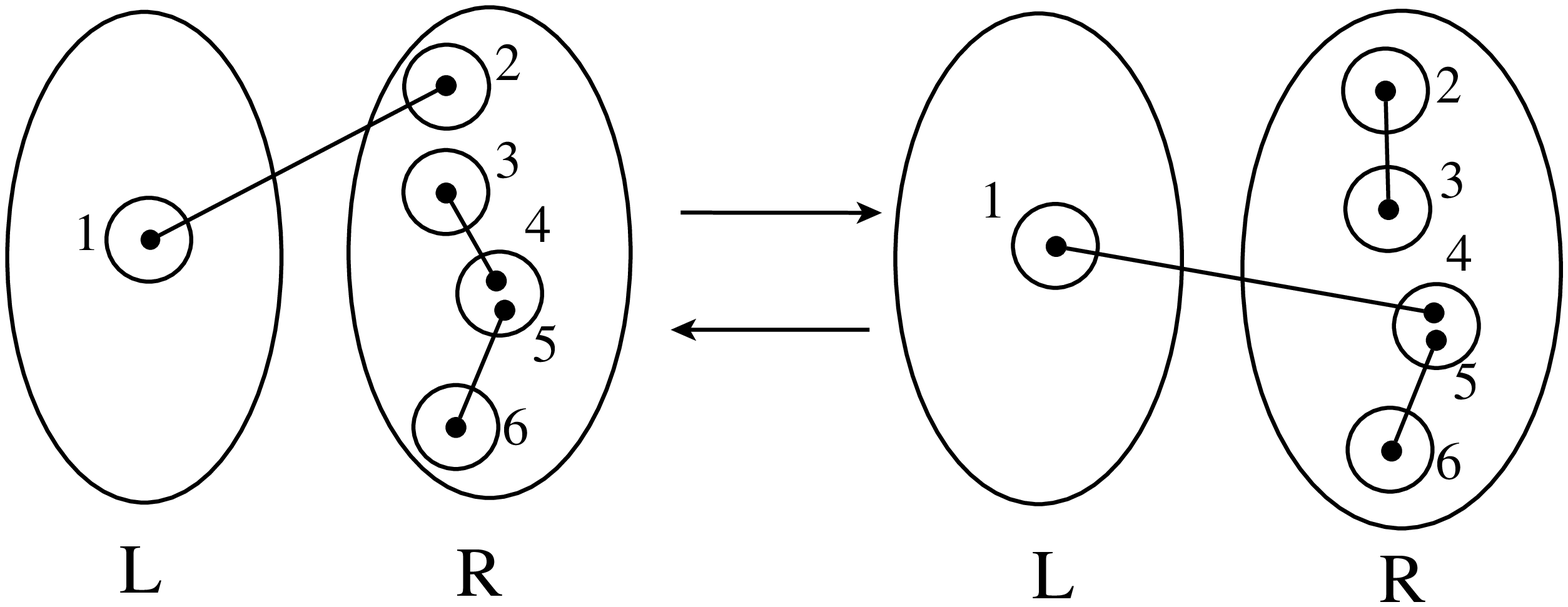}}}

\vskip .5cm \smallskip} \caption{\it $S_1$-switching}

\lab{f:swt4}

\end{figure}

\item{(b)} {\em $S_2$ switching}: Take a pure pair $\{5,6\}$ and a simple
directed $2$-path $((1,2),(3,4))$ such that the six points are
located in five distinct vertices shown as in Figure~\ref{f:swt5}.
Replace these three pairs by $\{1,2\}$, $\{3,5\}$ and $\{4,6\}$.
\end{description}

\begin{figure}[htb]
\vbox{\vskip .8cm
 \hbox{\centerline{\includegraphics[width=9cm]{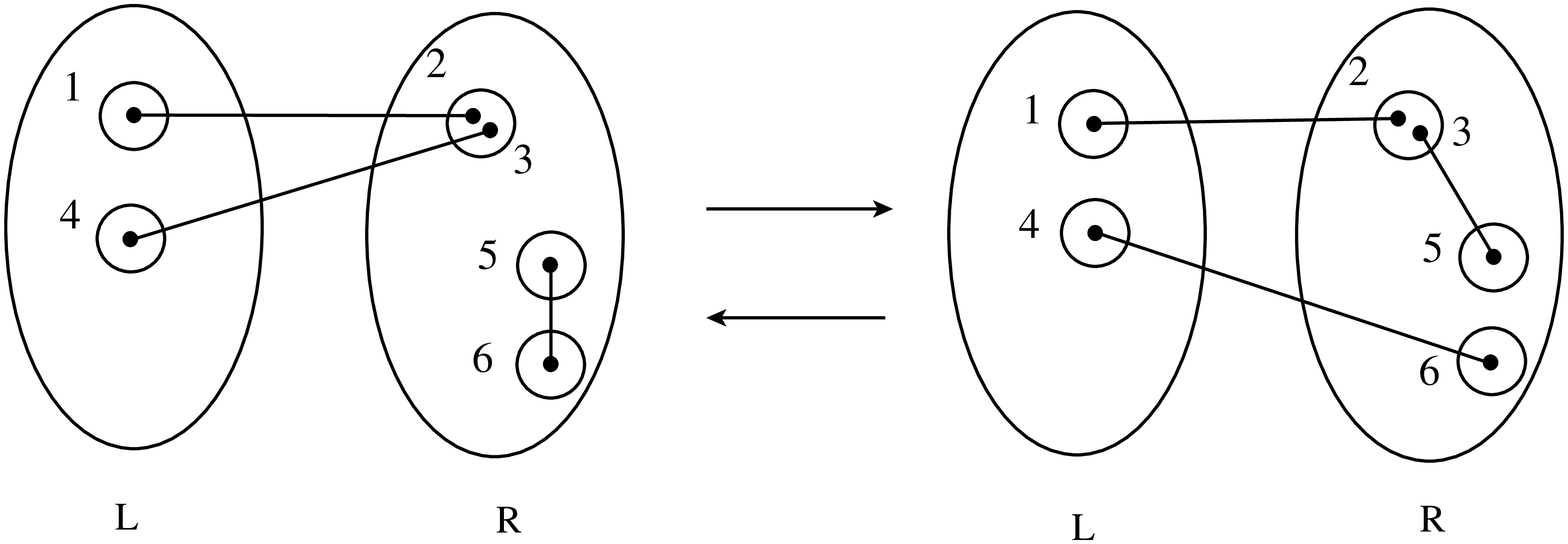}}}

\vskip .5cm \smallskip} \caption{\it  $S_2$-switching}

\lab{f:swt5}

\end{figure}
The inverse switchings are defined in the obvious way.
\begin{lemma}\lab{l:2Paths}
Given $l_0$, $l_1$ and $l_2$, let $\ell=l_0+l_1+l_2$.  We have
\begin{eqnarray*}
(i):&&\hspace{-0.4cm}\mbox{if}\  M_1(L)\le
M/4 \ \mbox{and}\ M_2(R)\ge 1,\\
&&\hspace{-0.5cm}a_1(l_0,l_1,l_2)=\frac{(M_1(R)-M_1(L))^2M_2(R)}{M_1(R)^2}(1+O(\d
^2/t+\ell/t+(\ell\d+l_0\d^2)/M_2(R)));\\
(ii):&&\hspace{-0.4cm}\mbox{if}\ M_1(L)>
M/4 \ \mbox{and}\ M_2(R)\ge 1,\\
&&\hspace{-0.5cm}a_3(l_0,l_1,l_2)=\frac{M_1(L)^2M_2(R)}{M_1(R)^2}(1+O(\d
^2/M_1(L)+\ell/M_1(L)+(\ell\d+l_0\d^2)/M_2(R))).
\end{eqnarray*}
\end{lemma}

\proof Let $a_i=a_i(l_0,l_1,l_2)$ for $i=1,2,3$. We use the
$S_1$-switching to compute the ratio $a_1/a_2$ and the
$S_2$-switching to compute the ratio $a_3/a_2$. We count the ordered
pairs of pairings $(\PA,\PA')$ such that both $\PA$ and $\PA'$ are
from $\mathC_{l_0,l_1,l_2}$, and $\PA'$ is obtained from $\PA$ by
applying an $S_1$-switching to $\PA$ without any creation or
destruction of loops or double pairs. Let $N_1$ denote the number of
such ordered pairs of pairings.

We first prove part (i). Assume $M_1(L)\le M/4$. For any directed
2-path of type 1 in $\mathC_{l_0,l_1,l_2}$, the number of
$S_1$-switching operations that can be applied to it is
\begin{equation}
A_1M_1(L)+O(A_1\d^2+A_1l_1)=A_1M_1(L)\Big(1+O(\d^2/M_1(L)+l_1/M_1(L))\Big).\lab{eq:s1}
\end{equation}
For any directed 2-path of type 2 in $\mathC_{l_0,l_1,l_2}$, the
number of inverse $S_1$-switching operations that can be applied to
it is
\begin{equation}
A_2\cdot 2t+O(A_2\d^2+A_2(l_0+l_2))=A_2\cdot
2t\Big(1+O(\d^2/t+(l_0+l_2)/t)\Big).\lab{eq:s2}
\end{equation}
The total number of $S_1$-switching operations that can be applied
to pairings in $\mathC_{l_0,l_1,l_2}$ is
$$
\sum_{\PA\in\mathC_{l_0,l_1,l_2}}A_1(\PA)M_1(L)\Big(1+O((\d^2+l_1)/M_1(L))\Big)=a_1M_1(L)\Big(1+O((\d^2+\ell)/M_1(L))\Big)|\mathC_{l_0,l_1,l_2}|,
$$
and the total number of inverse $S_1$-switching operations that can
be applied to pairings in $\mathC_{l_0,l_1,l_2}$ is
$$
\sum_{\PA\in\mathC_{l_0,l_1,l_2}}A_2(\PA)\cdot
2t\Big(1+O(\d^2/t+(l_0+l_2)/t)\Big)=a_2\cdot
2t\Big(1+O(\d^2/t+\ell/t)\Big)|\mathC_{l_0,l_1,l_2}|.
$$
These two numbers are both equal to $N_1$. Hence

\begin{eqnarray}
\frac{a_2}{a_1}&=&\frac{M_1(L)}{2t}(1+O(\d
^2/t+\d^2/M_1(L)+\ell/M_1(L)+\ell/t)). \lab{eq:ba}
\end{eqnarray}
Similarly, by the $S_2$-switching and its inverse we get
\begin{eqnarray}
\frac{a_3}{a_2}&=&\frac{M_1(L)}{2t}(1+O(\d ^2/t+\d
^2/M_1(L)+\ell/M_1(L)+\ell/t)).\lab{eq:bc}
\end{eqnarray}
Then~\eqn{eq:ba} gives
$$
\frac{a_2}{a_1}=\frac{M_1(L)}{2t}\Big(1+O((\d
^2+\ell)/t)\Big)+O((\d^2+\ell)/t),
$$
and~\eqn{eq:bc} gives
$$
\frac{a_3}{a_2} =\frac{M_1(L)}{2t}\Big(1+O((\d
^2+\ell)/t)\Big)+O((\d^2+\ell)/t).
$$
 Hence
\begin{eqnarray*}
a_2&=&a_1\left(\frac{M_1(L)}{2t}\Big(1+O((\d
^2+\ell)/t)\Big)+O((\d^2+\ell)/t)\right)\\
a_3&=&a_1\left(\frac{M_1(L)}{2t}\Big(1+O((\d
^2+\ell)/t)\Big)+O((\d^2+\ell)/t)\right)^2.
\end{eqnarray*}
Since $M_1(L)\le M/4$, we have $M_1(L)/t\le 1$ and so
$$
a_3=a_1\left(\left(\frac{M_1(L)}{2t}\right)^2\Big(1+O((\d
^2+\ell)/t)\Big)+O((\d^2+\ell)/t)\right).
$$
Hence
\begin{eqnarray}
a_1+2a_2+a_3\hspace{-0.2cm}&=&\hspace{-0.2cm}a_1\left(1+\left(2\frac{M_1(L)}{2t}+\left(\frac{M_1(L)}{2t}\right)^2\right)\Big(1+O((\d
^2+\ell)/t)\Big)+O((\d
^2+\ell)/t)\right)\nonumber\\
&=&\hspace{-0.2cm}a_1\left(\left(1+\frac{M_1(L)}{2t}\right)^2\Big(1+O((\d
^2+\ell)/t)\Big)+O((\d
^2+\ell)/t)\right)\nonumber\\
&=&\hspace{-0.2cm}a_1\left(1+\frac{M_1(L)}{2t}\right)^2(1+O((\d
^2+\ell)/t)).\lab{eq:ASum}
\end{eqnarray}

For any pairing $\PA\in\mathC_{l_0,l_1,l_2}$, the number of simple
directed $2$-paths in $\PA$ is $\sum_{v\in L\cup
R}d(v)(d(v)-1)-O(\ell\d+l_0\d^2)$, since the number of non-simple
directed 2-path is bounded by
$O(l_0\d^2+l_1\d+l_2\d)=O(\ell\d+l_0\d^2)$. On the other hand, the
number of simple directed $2$-paths in $\PA$ is $A_1+2A_2+A_3+A_4$,
since $2A_2$ counts the number of directed 2-paths of type 2 and the
opposite direction. Then
$$
A_1+2A_2+A_3+M_2(L)-O(l_1\d)=\sum_{v\in L\cup
R}d(v)(d(v)-1)-O(\ell\d+l_0\d^2).
$$
Thus,
\begin{equation}
A_1+2A_2+A_3=M_2(R)+O(\ell \d+l_0\d^2)=M_2(R)(1+O((\ell
\d+l_0\d^2)/M_2(R))).\lab{2path}
\end{equation}
 Combining this with~\eqn{eq:ASum}, we have
$$
a_1=\frac{(M_1(R)-M_1(L))^2M_2(R)}{M_1(R)^2}(1+O(\d
^2/t+\ell/t+(\ell\d+l_0\d^2)/M_2(R))),
$$
which proves part (i).

Next we show part (ii). Assume $M_1(L)>M/4$. We observe
that~\eqn{eq:ba} also gives
$$
\frac{a_1}{a_2}=\frac{2t}{M_1(L)}\Big(1+O((\d
^2+\ell)/M_1(L))\Big)+O((\d^2+\ell )/M_1(L)),
$$
and~\eqn{eq:bc} gives
$$
\frac{a_2}{a_3}=\frac{2t}{M_1(L)}\Big(1+O((\d
^2+\ell)/M_1(L))\Big)+O((\d^2+\ell )/M_1(L)).
$$
Thus,
\begin{eqnarray*}
a_2&=&a_3\left(\frac{2t}{M_1(L)}\Big(1+O((\d
^2+\ell)/M_1(L))\Big)+O((\d^2+\ell )/M_1(L))\right)\\
a_1&=&a_3\left(\frac{2t}{M_1(L)}\Big(1+O((\d
^2+\ell)/M_1(L))\Big)+O((\d^2+\ell )/M_1(L))\right).
\end{eqnarray*}
Since $M_1(L)\ge 1/4$, we have $t/M_1(L)<1$ and so
\begin{eqnarray*}
a_1+2a_2+a_3&=&a_3\left(1+\frac{2t}{M_1(L)}\right)^2\Big(1+O((\d
^2+\ell)/M_1(L))\Big)\\
&=& M_2(R)(1+O((\ell\d+l_0\d^2)/M_2(R))).
\end{eqnarray*}
Hence
$$
a_3=\frac{M_1(L)^2M_2(R)}{M_1(R)^2}(1+O(\d
^2/M_1(L)+\ell/M_1(L)+(\ell\d+l_0\d^2)/M_2(R))).
$$
This proves part (ii) of the lemma. \qed

\ss

 Recall the definition of $b_i(l_0,l_1,l_2)$ above
Lemma~\ref{l:double2}.  We next estimate these using simple
modifications of the $S_i$-switchings for $i=1,3$. (Note: in this lemma, our abbreviation $b_i$ contains no shift of index, whilst it did in Lemma~\ref{l:double2}.)
\begin{lemma}\lab{l:secondMoment}
 For $i=1,3$, let $a_i=a_i(l_0,l_1,l_2)$ and
$b_i=b_i(l_0,l_1,l_2)$,
  and let
$\ell=l_0+l_1+l_2$. Assume $M_2(R)\ge 1$. Then
\begin{eqnarray*}
(i): &&\mbox{if}\ M_1(L)\le M/4,\ \  b_1=a_1^2(1+O(\d^2/t+\ell/t+(\ell\d+l_0\d^2+\d^3)/M_2(R)));\\
(ii):&& \mbox{if}\ M_1(L)>M/4,\\
&&
b_3=a_3^2(1+O(\d^2/M_1(L)+\ell/M_1(L)+(\ell\d+l_0\d^2+\d^3)/M_2(R))).
\end{eqnarray*}
\end{lemma}

\proof For $1\le i\le 5$, let $X_i(\PA)$ denote the number of
ordered pairs of vertex disjoint simple $2$-paths in $\PA$ where the
first path has type $j_i$ and the second has type $h_i$, with
$(j_1,h_1)=(1,1)$,  $(j_2,h_2)=(3,3)$,  $(j_3,h_3)=(1,2)$,
$(j_4,h_4)=(1,3)$, and  $(j_3,h_3)=(2,3)$.

 The $S_3$-switching, as illustrated in Figure~\ref{f:swt6},  is a slight
modification of the $S_1$-switching. To apply it, we need to choose
a mixed pair and two simple $2$-paths of type 1 such that they are
pairwise disjoint. To apply its inverse, we need to choose a pure
pair and two simple $2$-paths of type 2 and 1 respectively such that
they are pairwise disjoint. Compared with the $S_1$-switching, the
$S_3$-switching   requires  an additional simple directed $2$-path
of type $1$. However, the pairs in the extra $2$-path remain after
the $S_3$-switching is applied since the mixed pair and the other
simple directed $2$-path under consideration are vertex-disjoint
from the additional directed $2$-path. The $S_4$-switching, as
 illustrated in Figure~\ref{f:swt7}, is a similar modification of the $S_2$-switching.

We will first estimate $\ex(X_i(\PA)\mid
\PA\in\mathC_{l_0,l_1,l_2})$ for $i\in [5]$ and then use this to
estimate $b_1$ and $b_3$. Following the analogous argument as in
Lemma~\ref{l:2Paths}, we can estimate the ratio  $\ex(X_3(\PA)\mid
\PA\in\mathC_{l_0,l_1,l_2})/\ex(X_1(\PA)\mid
\PA\in\mathC_{l_0,l_1,l_2})$ by counting the ordered pairs of
pairings $(\PA,\PA')$ such that $\PA,\PA'\in\mathC_{l_0,l_1,l_2}$
and $\PA'$ is obtained by applying an $S_3$-operation to $\PA$
without any creation or destruction of loops or double pairs.
 The the number of such $S_3$-switching
operations that can be applied to $\PA$ is
$X_1M_1(L)+O(X_1\d^2+X_1l_1)$. The number of such inverse
$S_3$-operations that can be applied to $\PA$ is
$2tX_3+O(X_3\d^2+X_3(l_0+l_2))$.  So the ratio  $\ex(X_3(\PA)\mid
\PA\in\mathC_{l_0,l_1,l_2})/\ex(X_1(\PA)\mid
\PA\in\mathC_{l_0,l_1,l_2})$ equals exactly the right hand side
of~\eqn{eq:ba} and the ratio $\ex(X_4(\PA)\mid
\PA\in\mathC_{l_0,l_1,l_2})/\ex(X_3(\PA)\mid
\PA\in\mathC_{l_0,l_1,l_2})$ equals exactly the right hand side
of~\eqn{eq:bc}.

\begin{figure}[htb]
\vbox{\vskip .8cm
 \hbox{\centerline{\includegraphics[width=10cm]{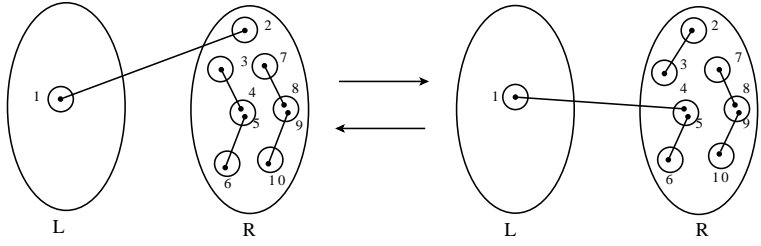}}}

\vskip .5cm \smallskip} \caption{$S_3$-switching}

\lab{f:swt6}

\end{figure}

\begin{figure}[htb]
\vbox{\vskip .8cm
 \hbox{\centerline{\includegraphics[width=10cm]{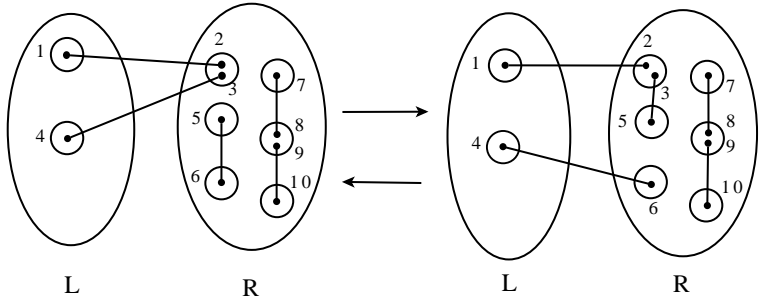}}}

\vskip .5cm \smallskip} \caption{$S_4$-switching}

\lab{f:swt7}

\end{figure}

%
On the other hand, by~\eqn{2path}, for any
$\PA\in\mathC_{l_0,l_1,l_2}$,
$A_1(\PA)+2A_2(\PA)+A_3(\PA)=M_2(R)(1+O((\ell\d+l_0\d^2)/M_2(R)))$.
Thus,
\begin{eqnarray}
&&\ex(A_1^2\mid \PA\in\mathC_{l_0,l_1,l_2})+2\ex(A_1A_2\mid \PA\in\mathC_{l_0,l_1,l_2})+\ex(A_1A_3\mid \PA\in\mathC_{l_0,l_1,l_2})\non\\
&&\hspace{.3cm}=\ex(A_1(A_1+2A_2+A_3)\mid \PA\in\mathC_{l_0,l_1,l_2})\non\\
&&\hspace{.3cm}=\ex(A\mid \PA\in\mathC_{l_0,l_1,l_2})M_2(R)(1+O((\ell\d+l_0\d^2)/M_2(R)))\non\\
&&\hspace{.3cm}=a_1(l_0,l_1,l_2)M_2(R)(1+O((\ell\d+l_0\d^2)/M_2(R))).\lab{2path2}
\end{eqnarray}
We also have
\begin{equation} X_1=A_1^2+O(A_1\d^3),\ \ X_3=A_1A_2+O(A_1\d^3),\ \
X_4=A_1A_3+O(A_1\d^3),\lab{eq:disjoint}
\end{equation}
where the error terms in~\eqn{eq:disjoint} account for the number of
ordered pairs of simple $2$-directed paths that are not vertex
disjoint. Let $a_1=a_1(l_0,l_1,l_2)$. Taking the conditional
expectation on both sides of each equation in~\eqn{eq:disjoint}, we
obtain
\begin{eqnarray*}
\ex(X_1\mid \PA\in\mathC_{l_0,l_1,l_2})&=&\ex(A_1^2\mid
\PA\in\mathC_{l_0,l_1,l_2})+O(a_1\d^3),\\
\ex(X_3\mid \PA\in\mathC_{l_0,l_1,l_2})&=&\ex(A_1A_2\mid
\PA\in\mathC_{l_0,l_1,l_2})+O(a_1\d^3),\\
\ex(X_4\mid \PA\in\mathC_{l_0,l_1,l_2})&=&\ex(A_1A_3\mid
\PA\in\mathC_{l_0,l_1,l_2})+O(a_1\d^3).
\end{eqnarray*}
Combining this with~\eqn{2path2} we have
\begin{eqnarray*}
&&\ex(X_1\mid \PA\in\mathC_{l_0,l_1,l_2})+2\ex(X_3\mid
\PA\in\mathC_{l_0,l_1,l_2})+\ex(X_4\mid
\PA\in\mathC_{l_0,l_1,l_2})\\
&&=a_1M_2(R)(1+O((\ell\d+l_0\d^2+\d^3)/M_2(R))).
\end{eqnarray*}
So part (i) follows from an argument similar to that used for
Lemma~\ref{l:2Paths} and~\eqn{eq:disjoint}.
Similarly, by analysing two switching operations similar to those of
$S_3$-switching and $S_4$-switching, except that the extra $2$-path
is of type 3, we can estimate the ratio
$\ex(X_5\mid\PA\in\mathC_{l_0,l_1,l_2})/\ex(X_4\mid\PA\in\mathC_{l_0,l_1,l_2})$
and
$\ex(X_2\mid\PA\in\mathC_{l_0,l_1,l_2})/\ex(X_4\mid\PA\in\mathC_{l_0,l_1,l_2})$.
By the fact that
$$
X_5=A_2A_3+O(A_3\d^3),\ \ X_4=A_1A_3+O(A_3\d^3),\ \
X_2=A_3^2+O(A_3\d^3),
$$
and
\begin{eqnarray*}
&&\ex(X_2\mid \PA\in\mathC_{l_0,l_1,l_2})+2\ex(X_5\mid
\PA\in\mathC_{l_0,l_1,l_2})+\ex(X_4\mid
\PA\in\mathC_{l_0,l_1,l_2})\\
&&=a_3(l_0,l_1,l_2)M_2(R)(1+O((\ell\d+l_0\d^2+\d^3)/M_2(R))),
\end{eqnarray*}
  together
with Lemma~\ref{l:2Paths}(ii), part (ii) follows from an argument
similar to that in part (i) and the proof of
Lemma~\ref{l:2Paths}(ii).  \qed

\newpage
\section{Synthesis}
\lab{s:synthesis} We are now ready to substitute the values of the
variables $a_i$ and $b_i$ determined in
Section~\ref{s:moreswitchings} in the ratios determined in
Section~\ref{s:basics}, and from there to prove the main theorem.
The reader should not be surprised at how the separate cases combine
to give the same resulting formulae with the desired error terms;
the definitions of the cases and the choices of switchings for each
case were carefully designed to achieve this.
\begin{lemma}\lab{l:ratio}
Assume $\d^4=o(M)$. Let  $\alpha_0=((l_1+l_2)\d+l_0\d^2)/M_2(R)$,
$\alpha_1=((l_1+l_2)\d)/M_2(R)$ and $\alpha_2=(l_2\d+\d^3)/M_2(R)$.
Assume $M_2(R)/\d^3$ is sufficiently large. Then
\begin{eqnarray*}
(i)\ \frac{|\mathC_{l_0,l_1,l_2}|}{|\mathC_{l_0-1,l_1,l_2}|}\hspace{-0.3cm}&=&\hspace{-0.3cm}\frac{\mu_0}{l_0}\left(1+O\left(\frac{\d^2+l_0+l_2}{t}+\frac{l_1}{M}\right)\right)(1+O(\alpha_0)),\  l_0\ge 1;\\
(ii)\ \frac{|\mathC_{0,l_1,l_2}|}{|\mathC_{0,l_1-1,l_2}|}\hspace{-0.3cm}&=&\hspace{-0.3cm}\frac{\mu_1}{l_1}\left(1+O\left(\frac{\d^3+l_1\d}{M_2(L)}+\frac{\d^2+l_1+l_2}{M}\right)\right)(1+O(\alpha_1)),\  l_1\ge 1;\\
(iii)\
\frac{|\mathC_{0,0,l_2}|}{|\mathC_{0,0,l_2-1}|}\hspace{-0.3cm}&=&\hspace{-0.3cm}\frac{\mu_2}{l_2}\left(1+O\left(\frac{\d^3}{M_2(R)}+\frac{\d^2}{M}+\frac{l_2}{t}\right)\right)(1+O(\alpha_2)),\
l_2\ge 1.
\end{eqnarray*}

\end{lemma}
\proof Let $\delta=M_1(L)/M$, so that $0\le \delta\le 1/2$ since $M_1(R)\ge M_1(L)$.
\smallskip

\no {\em Case 1}: $ \delta\le 1/4$.

Here $t$, which was defined as $(M_1(R)-M_1(L))/2$, is $\Theta(M)$.
By part (i) of Lemmas~\ref{l:loop}--\ref{l:double2}
and~\ref{l:2Paths}--\ref{l:secondMoment}, and
recalling~\eqn{eq:mu0}--\eqn{eq:mu2}, we obtain the following, with
some of the bounds on error terms explained below.
\begin{eqnarray*}
\frac{|\mathC_{l_0,l_1,l_2}|}{|\mathC_{l_0-1,l_1,l_2}|}&=&\frac{\mu_0}{l_0}(1+O(\d ^2/M+(l_0+l_2)/M))(1+O(\alpha_0)),\\
\frac{|\mathC_{0,l_1,l_2}|}{|\mathC_{0,l_1-1,l_2}|}&=&
\frac{\mu_1}{l_1}(1+O((\d
^3+l_1\d)/M_2(L)+(\d^2+l_2)/M))(1+O(\alpha_1)),\\
\frac{|\mathC_{0,0,l_2}|}{|\mathC_{0,0,l_2-1}|}&=&
\frac{\mu_2}{l_2}(1+O((\d ^2+l_2)/M+\d^3a_1/b_1))(1+O(\alpha_3)).
\end{eqnarray*}
For the second equations, note that error terms involving $l_0$ do
not appear since $l_0=0$, and similarly $l_0=l_1=0$ for the third
equation.
\smallskip

\no {\em Case 2}: $1/4<\delta\le 1/2$.

Here
 $M_1(L)=\Theta(M)$. By part (ii) of Lemmas~\ref{l:loop}--\ref{l:double2} and~\ref{l:2Paths}--\ref{l:secondMoment},
  we obtain the following, with some error terms explained below.
\begin{eqnarray*}
\frac{|\mathC_{l_0,l_1,l_2}|}{|\mathC_{l_0-1,l_1,l_2}|}
&=&
\frac{\mu_0}{l_0}(1+O((\d^2+l_1)/M+(\d^2+l_0+l_2)/t))(1+O(\alpha_0)),\\
\frac{|\mathC_{0,l_1,l_2}|}{|\mathC_{0,l_1-1,l_2}|}&=&
\frac{\mu_1}{l_1}(1+O((\d^2+l_1)/M)+(\d^3+l_1\d)/M_2(L))(1+O(\alpha_1)),\\
\frac{|\mathC_{0,0,l_2}|}{|\mathC_{0,0,l_2-1}|}&=&
\frac{\mu_2}{l_2}(1+O(\d ^3a_3/b_3)+\d^2/M+l_2/t)(1+O(\alpha_2)).
\end{eqnarray*}
To obtain the second of these equations, note that $l_1/M_1(L)=O(l_1/M)=O(l_1\d/M_2(R))=O( \alpha_1)$.

Parts (i) and (ii) follow by combining the two cases. To complete
the proof of part (iii), we show that $a_1/b_1=O(M_2(R)^{-1})$ when
$M_1(L)\le M/4$ and $a_3/b_3=O(M_2(R)^{-1})$ when $M_1(L)>M/4$.

First consider $M_1(L)\le M/4$. Considering $a_1/b_1$, we have the following two
cases.
\smallskip

\no{\em Case 1}: $M_2(R)\le \zeta_2(\d^5+\d^3\ln^2 M)$. Then
$k_2=\d^2+2$ according to its redefinition  after
Lemma~\ref{c:expectation}. Since $M_2(R)/\d^3$ can be assumed
arbitrarily large by the present lemma's assumption, the error terms
$ l_2\d/M_2(R)$ and $\d^3/M_2(R)$ in Lemmas~\ref{l:2Paths}(i)
and~\ref{l:secondMoment}(i) can be taken arbitrarily small. It
follows that $a_1=\Omega(M_2(R))$ and $b_1=\Theta(a_1^2)$, and so
$a_1/b_1=O(M_2(R)^{-1})$.
\smallskip

\no{\em Case 2}: $M_2(R)/(\d^5+\d^3\ln^2 M)>\zeta_2$, which can at
this point be taken arbitrarily large. Then for any $l_2\le
k_2=O(\d^2)$, as defined in~\eqn{kdefs}, the error terms in
Lemmas~\ref{l:2Paths}(i) and~\ref{l:secondMoment}(i) can be made
arbitrarily small. Thus  $a_1=\Omega(M_2(R))$,  $b_1=\Theta(a_1^2)$,
and $a_1/b_1=O(M_2(R)^{-1})$.

On the other hand, assuming $M_1(L)>M/4$, a similar argument shows
that $a_3/b_3=O(M_2(R)^{-1})$. \qed \ss

Recall that $\pr( {\bf d})$ denotes the probability that a random
pairing $\PA\in \mathcal{M}(L,R,{\bf d})$ corresponds to a simple
B-graph.
\ss

\no {\bf Proof of Theorem~\ref{t:graphCount}.\ } Recall that   $\pr( {\bf d})$ denotes the probability that a random
pairing $\PA\in \mathcal{M}(L,R,{\bf d})$ corresponds to a simple
B-graph, and
$\mathscrC(m)$ denotes the number $m!/\big((m/2)!2^{m/2}\big)$ of pairings of $m$ points. The
total number of pairings in $\mathcal{M}(L,R,{\bf d})$ is thus
$
[M_1(R)]_{M_1(L)}\mathscrC(M_1(R)-M_1(L))$.
Since each simple B-graph corresponds to $\prod_{i=1}^n d_i$
pairings in $\mathcal{M}(L,R,{\bf d})$, we have
$$
g(L,R,{\bf d})= \frac{M_1(R)!\pr({\bf
d})}{2^{(M_1(R)-M_1(L))/2}((M_1(R)-M_1(L))/2)!\prod_{i=1}^n d_i!},
$$
and it only remains to show that
$\pr({\bf d})=e^{-\mu_0-\mu_1-\mu_2}(1+O(\d^4/M)).
$

If $M_2(R)=O(\d^3)$, we have $\mu_i=O(\d^4/M)$ for $i=0,1,2$. Then
by Corollary~\ref{c2:expectation} and the first moment principle,
$\pr({\bf d})=1-O(\d^4/M)$ and we are done. So we may assume
\bel{M2big} M_2(R)/\d^3>C \ee for any arbitrarily large $C$. (Note
we   assume throughout that  $\d>0$ since otherwise there is nothing
to prove.)  By Corollary~\ref{c:probability}, it is enough to show
\bel{dosyc} \sum_{l_2=0}^{k_2}
\sum_{l_1=0}^{k_1}\sum_{l_0=0}^{k_0}|\mathC_{l_0,l_1,l_2}|=|\mathC_{0,0,0}|e^{\mu_0+\mu_1+\mu_2}(1+O(\d^4/M)).
\ee
Iterating the ratio in Lemma~\ref{l:ratio}(i), for any fixed $l_0\le k_0$, $l_1\le k_1$ and $l_2\le k_2$, we get
 $$
\frac{|\mathC_{l_0,l_1,l_2}|}{|\mathC_{0,l_1,l_2}|}=\frac{\mu_0^{l_0}}{l_0!}\left(1+O(\d^2/t+(l_0+l_2)/t+l_1/M)\right)^{l_0}(1+O(\alpha_0))^{l_0}
 $$
where $\alpha_0$ is as defined in that lemma.

First we sum over $l_0$. Here we assume $t\ge 1$,  since otherwise  $B_0=0$, which will trivially give the desired conclusion. Recalling
the definition~\eqn{kdefs} of $k_i$ and its redefinition after
Corollary~\ref{c:expectation}, we have $k_0=O(\d+\ln M)$ and for
$i=1,2$, $k_i=O(\d^2+\ln M)$. Consider the following two cases,
recalling $t$ from~\eqn{tdef}. \ss

 \no {\em Case 1}:  $M_2(R)\le \zeta_0(\d^5+\d^3\ln^2 M)$ or $2t\le  \zeta_1(\d^4+\d^2\ln^2 M)$.
\ss

\no Here, by the redefinition of $k_i$, we have $k_0=O(\d)$ and
$k_2=O(\d^2)$, so $ \alpha_0 =  (O(\d^3/M_2(R)))$. Recalling also the definition~\eqn{eq:mu0} of $\mu_0$ as $tM_2(R)/M_1(R)^2$, and noting $M_1(R)=\Omega(M)$ and $M_2(R)=O(\d  M)$, we have from Lemma~\ref{l:ratio}(i) that for $1\le l_0\le k_0$ and all relevant $l_1$ and $l_2$,
\[
 \frac{|\mathC_{l_0,l_1,l_2}|}{|\mathC_{l_0-1,l_1,l_2}|}
 =
   \frac{1}{l_0}\big(\mu_0/l_0 +O(\d^3/M) +O(\d l_1/M)).
\]
  Hence (bounding $\d^3$ by $\d^4$  for consistency with   the later argument),
\begin{eqnarray*}
\sum_{l_0=0}^{k_0} \frac{|\mathC_{l_0,l_1,l_2}|}{|\mathC_{0,l_1,l_2}|}
&=& \sum_{l_0=0}^{k_0}\frac{(\mu_0+O(\d^4/M)+O(\d l_1/M))^{l_0}}{l_0!}\\
&=&  \exp\big(\mu_0+O(\d^4/M+\d l_1/M)\big)+O\big((\d^4+\d l_1)/M\big)
\end{eqnarray*}
using
$$
\sum_{l_0=k_0+1}^{\infty}\frac{ (\mu_0+x )^{l_0}}{l_0!}=
\sum_{l_0=k_0+1}^{\infty}\frac{ \big(O(\mu_0)\big)^{l_0}+\big(O(x ))^{l_0}}{l_0!}=
O(\mu_0^{k_0}/k_0!+x)
$$
for $x=o(1)$, and noting that $\mu_0=O(\d^5/M)$ in this case, which is $o(\d)$ and hence less than $\d/2$ for large $M$. (In particular,  $\mu_0$ tends to 0 quickly unless $\d$ is large.) Hence
\[
\sum_{l_0=0}^{k_0} \frac{|\mathC_{l_0,l_1,l_2}|}{|\mathC_{0,l_1,l_2}|}
=\exp(\mu_0)\big(1+O(\d^4/M+\d l_1/M)\big).
\]

\no {\em Case 2}:
$M_2(R)>\zeta_0(\d^5+\d^3\ln^2 M)$ and $2t>\zeta_1(\d^4+\d^2\ln^2
M)$.

\no Here $k_0=O(\ln M+\d)$, $k_i=O(\ln M+\d^2)$ for $i=1,2$.  Note that $\d^3\ln M\le \d^4+\d^2\ln^2 M =O(t)$, and from here we see that $k_0\d^2/t=O(1)$. Similarly, $k_0k_2  = O(\ln^2M +\d^3)=O(t)$. In this way, we find that
$l_0(\d^2/t+(l_0+l_2)/t+l_1/M+\alpha_0)=O(1)$ provided $l_i\le k_i$ for $i=0,1,2$.
So, from Lemma~\ref{l:ratio}(i),
\begin{eqnarray*}
 \sum_{l_0=0}^{k_0}\frac{|\mathC_{l_0,l_1,l_2}|}{|\mathC_{0,l_1,l_2}|}&=&\sum_{l_0=0}^{k_0}\frac{\mu_0^{l_0}\exp(O(l_0(\d^2/t+(l_0+l_2)/t+l_1/M+\alpha_0)))}{l_0!}\\
& =&\sum_{l_0=0}^{k_0}\frac{\mu_0^{l_0}}{l_0!}+O\left(\sum_{l_0=0}^{k_0}\frac{\mu_0^{l_0}}{l_0!}l_0\left(\frac{\d^2+l_2}{t}+\frac{l_1}{M}+\frac{(l_1+l_2)\d}{M_2(R)}\right)\right)\\
&&\hspace{0.6cm}+O\left(\sum_{l_0=0}^{k_0}\frac{\mu_0^{l_0}}{l_0!}l_0^2\left(\frac{1}{t}+\frac{\d^2}{M_2(R)}\right)\right).
\end{eqnarray*}
Note also that $k_0\ge 8\eta(R)\ge 16\mu_0$, and $k_0\ge \ln M$. So
$$
\sum_{l_0=k_0+1}^{\infty}\frac{\mu_0^{l_0}}{l_0!}=O((k_0/16)^{k_0}{k_0!}=O\big((e/16)^{k_0}\big)=o(M^{-1}).
$$
Also, of course, $\sum_{l_0=0}^{k_0}(\mu_0^{l_0}/l_0!)l_0\le \mu_0e^{\mu_0}$ and $\sum_{l_0=0}^{k_0}(\mu_0^{l_0}/l_0!)l_0^2\le (\mu_0^2+\mu_0)e^{\mu_0}$. So we have
\begin{eqnarray*}
 \sum_{l_0=0}^{k_0}\frac{|\mathC_{l_0,l_1,l_2}|}{|\mathC_{0,l_1,l_2}|}
&=&e^{\mu_0}-O(M^{-1})+O\left(e^{\mu_0}\mu_0\left(\frac{\d^2+l_2}{t}+\frac{l_1}{M}+\frac{(l_1+l_2)\d+\d^3}{M_2(R)}\right)\right)\\
&&\hspace{0.6cm}+O\left(e^{\mu_0}(\mu_0^2+\mu_0)\left(\frac{1}{t}+\frac{\d^2}{M_2(R)}\right)\right).
\end{eqnarray*}
Now using
\bean
\mu_0/t&=&M_2(R)/M_1(R)^2=O(\d/M),\\
\mu_0^2/t&=&O(M_2(R)^2t/M_1(R)^4)=O(\d^2/M_1),\\
\mu_0  &=&O(M_2(R)/M),\\
\mu_0  &=&O(\d),
\eean
we obtain
\[
 \sum_{l_0=0}^{k_0}\frac{|\mathC_{l_0,l_1,l_2}|}{|\mathC_{0,l_1,l_2}|}
 = e^{\mu_0}\left(1+O\left(\frac{(l_1+l_2)\d}{M}+\frac{\d^3}{M}\right)\right).
\]
Combining the two cases, we have (for $l_1$ and $l_2$ in the appropriate range)
$$
\sum_{l_0=0}^{k_0}|\mathC_{l_0,l_1,l_2}|=|\mathC_{0,l_1,l_2}|\exp(\mu_0)\left(1+O\left(\frac{(l_1+l_2)\d}{M}+\frac{\d^4}{M}\right)\right).
$$

We will next sum this expression over $l_1$.
 By
Lemma~\ref{l:ratio}(ii), for any fixed $l_1\le k_1$ and $l_2\le
k_2$,
$$
\frac{|\mathC_{0,l_1,l_2}|}{|\mathC_{0,0,l_2}|}=\frac{\mu_1^{l_1}}{l_1!}\left(1+O\left(\frac{\d^3+l_1\d}{M_2(L)}+\frac{\d^2}{M}+\frac{l_1+l_2}{M}\right)\right)^{l_1}(1+O(\alpha_1))^{l_1}
$$
where   $\alpha_1=(l_1+l_2)\d/M_2(R)$.
\medskip

\no {\em Case 1}: $M_2(R)\le \zeta_0(\d^5+\d^3\ln^2 M)$ or
$M_2(L)\le\zeta_1(\d^5+\d^3\ln^2 M)$.
Then $k_1=\d^2+2$, and summing over
$0\le l_1\le k_1$ we obtain
$$
\sum_{l_1=0}^{k_1}\sum_{l_0=0}^{k_0}|\mathC_{l_0,l_1,l_2}|=\exp(\mu_0+\mu_1)|\mathC_{0,0,l_2}|\left(1+O\left(\frac{l_2\d^2+\d^4}{M}\right)\right).
$$
\no{Case 2}: $M_2(R)>\zeta_0(\d^5+\d^3\ln^2 M)$ and
$M_2(L)>\zeta_1(\d^5+\d^3\ln^2 M)$. Then for any $l_1\le k_1$,
$l_2\le k_2$,
$$
l_1\left(\frac{\d^3+l_1\d}{M_2(L)}+\frac{\d^2}{M}+\frac{l_1+l_2}{M}+\alpha_1\right)
$$
is bounded. Estimating error terms similar to Case~2 of the earlier summation over $l_0$, we
obtain the same result as in Case~1.

 For summing over $l_2$, the argument is similar, and the final result is~\eqn{dosyc} as
 required.\qed\ss

\end{document}